\documentclass[11pt]{article}

\def\squarebox#1{\hbox to #1{\hfill\vbox to #1{\vfill}}}
\def\boxit#1{\vbox{\hrule\hbox{\vrule\kern6pt
          \vbox{\kern6pt#1\kern6pt}\kern6pt\vrule}\hrule}}

\usepackage{graphicx}
\RequirePackage{amsmath,amsthm,amsfonts,amssymb}
\usepackage{dsfont}
\usepackage{comment}
\usepackage{bm}
\usepackage{color,soul}
\usepackage{multirow}
\usepackage{natbib}
\usepackage{enumerate}
\usepackage{subcaption}
\usepackage{mathtools}

\RequirePackage[citecolor=blue]{hyperref}

\renewcommand{\baselinestretch} {1.5}
\makeatletter 
\def\singlespace{\def\baselinestretch{1}\@normalsize}

\date{}

\renewcommand{\theequation} {\arabic{section}.\arabic{equation}}
\newtheorem{theorem}{Theorem}
\newtheorem{lemma}{Lemma}%[section]
\newtheorem{remark}{Remark}%[section]
\newtheorem{corollary}{Corollary}%[section]

\newtheorem{example}{Example}
\newtheorem{assumption}{Assumption}

\theoremstyle{plain}
\newtheorem{definition}{Definition}

\def\R{{\mathbb R}}

\def\N{{\mathbb N}}

\def\E{{\mathbb E}}
\def\R{{\mathbb R}}

\def\L{{\mathcal L}}

\def\D{{\mathcal D}}
\def\P{{\mathbb P}}

\def\A{{\mathcal A}}
\def\W{{\bm W}}

\def\cG{{\mathcal G}}
\def\cE{{\mathcal E }}
\def\cV{{\mathcal V }}

\def\G{{\mathsf G}}

\def\u{{\bm{u}}}

\def\w{{\bm{w}}}
\def\uu{{\widehat{\bm{u}}}}

\def\bbA{{\mathbb A}}

\DeclareMathOperator*{\argmax}{arg\,max}

\renewcommand{\d}[1]{\ensuremath{\operatorname{d}\!{#1}}}

\begin{document}
\baselineskip=22pt
\title{\bf\Large Statistical inference for pairwise comparison models}
%\title{\bf\Large Inference in pairwise comparison models}
\author{
{Ruijian Han\thanks{Authorships are ordered alphabetically. } \textsuperscript{1}, Wenlu Tang$^*$\textsuperscript{2} and
Yiming Xu$^*$\textsuperscript{3}}\\
\small \textsuperscript{1} Department of Data Science and Artificial Intelligence, The Hong Kong Polytechnic University, \\ [-3mm] \small Hung Hom, Kowloon, Hong Kong SAR
\\[-3mm]
\small \textsuperscript{2}Department of Mathematical and Statistical Sciences, University of Alberta,   Edmonton, Alberta, Canada 
\\
[-3mm]
\small{\textsuperscript{3} Department of Mathematics, University of Kentucky, Lexington, Kentucky, USA} 
}
%\date{\today}
\maketitle

\begin{center}
\bf \Large Abstract
\end{center}
Pairwise comparison models have been widely used for utility evaluation and rank aggregation across various fields. The increasing scale of modern problems underscores the need to understand statistical inference in these models when the number of subjects diverges, a topic that is currently underexplored in the literature. To address this gap, this paper establishes a near-optimal asymptotic normality result for the maximum likelihood estimator in a broad class of pairwise comparison models. The key idea lies in identifying the Fisher information matrix as a weighted graph Laplacian, which can be studied via a meticulous spectral analysis. Our findings provide theoretical foundations for performing statistical inference in a wide range of pairwise comparison models beyond the Bradley--Terry model. Simulations utilizing synthetic data are conducted to validate the asymptotic normality result, followed by a hypothesis test using a tennis competition dataset.

		\vspace{0.1in}


\newpage

	\section{Introduction}

Pairwise comparison involves assessing subjects in pairs to establish their relative preferences, a practice pertinent to various applications such as econometrics \citep{eco1958,mcfadden1973conditional,10.1257/pandp.20221063,lentz2023anatomy}, sports analytics \citep{baker2014dynamic, bozoki2016application,angelini2022weighted,collingwood2022evaluating}, social science \citep{loewen2012testing,10.1111/rssa.12124,wapman2022quantifying}, and, more recently, human preference studies aimed at improving models in artificial intelligence \citep{Christiano2017, Rafailov2023, Sun2025}.  
 A prevalent approach for pairwise comparison modeling employs a latent score framework. Originated from the ideas in Thurstone \citep{thurstone1927method} and Zermelo \citep{MR1545015}, a mathematical model for pairwise comparison data analysis was formulated by Bradley and Terry \citep{MR0070925}. Since then, multiple generalizations have been developed, including ordinal models such as the Rao--Kupper model \citep{MR217963} and the Davidson model \citep{davidson1970extending}, which account for ties, the cumulative link model \citep{agresti1992analysis} that considers more refined ordinal scales, and cardinal models such as the paired cardinal model \citep{MR3504618}. 
 We recommend \cite{MR3012434} for a review of pairwise comparison modeling from a practical perspective. 

In light of the fast-growing scale of problems in the big-data era, a recent research area of pairwise comparison models focuses on understanding the asymptotics of estimating the latent score vector as the number of compared subjects approaches infinity.
Along this line of work, considerable attention has been devoted to the study of the Bradley--Terry (BT) model, which is popular due to its simple parametrization form and connections to other problems such as matrix balancing \citep{qu2023sinkhorn} and completion \citep{chen2023note}.
Notably, an efficient Minorize-Maximization (MM) algorithm for implementing the maximum likelihood estimator (MLE) was proposed in \cite{MR2051012} and accelerated in \cite{newman2023efficient} in the BT model, and its uniform consistency was established when the comparison graph is dense \citep{MR1724040, MR2987494}.
Subsequent works have established analogous consistency results for alternative estimators such as the spectral estimator and regularized MLE in the sparse settings under various error metrics \citep{chen2015spectral, MR3504618, MR3613103, MR3953449}. Recently, \cite{han2020asymptotic} and \cite{chen2022optimal} extended the uniform consistency of the MLE to the sparse setting, paralleled by a few other works with a slightly different focus (e.g., asymptotic efficiency and minimax rates); see \cite{hendrickx2020minimax, bong2022generalized}. It is noteworthy that \cite{han2020asymptotic} also contains an asymptotic normality result for dense graphs, which was further explored in the sparse setting \citep{MR4560195, gao2023uncertainty}.

Despite extensive research on the asymptotic properties of the BT model, a gap persists between theory and practice as many comparison models used in applications either have multiple scales or a continuous outcome. This suggests a practical need to extend this understanding to other widely used yet less explored models. The recent work \citep{han2022general} partially addresses this by introducing a general pairwise comparison model framework and proving a uniform consistency result for the MLE therein. However, an asymptotic normality result has not yet been established. Such a result is often crucial for practical inference and \emph{cannot} be deduced from the existing approaches to analyzing the BT model.

To take one step further in filling the gap, this paper aims to derive an asymptotic normality result for the MLE in the general pairwise comparison models. In particular, our contributions can be summarized as follows:
\begin{itemize}
	\item We establish an asymptotic normality result for the MLE in the general pairwise comparison models under a near-optimal sparsity condition on the comparison graphs. Additionally, we derive a non-asymptotic convergence rate for each individual subject under comparison. The crux of our analysis involves recognizing the Fisher information matrix as a weighted graph Laplacian, enabling a meticulous analysis based on spectral expansion. 
	Aside from enjoying the sufficient flexibility of model parametrization, our results additionally allow the incorporation of imbalanced data, which is new compared to many existing results in the field and may be of independent interest; see Table~\ref{tab:my_label} for a detailed comparison between our result and some existing results in the literature. 
	\item We verify that the MLE in many pairwise comparison models in the literature, including the BT model, the Thurstone--Mosteller model \citep{thurstone1927method, mosteller2006remarks}, the Davidson model, the Rao--Kupper model, the cumulative link model, and the paired cardinal model, is asymptotically normal under near-optimal sparsity conditions.
	We also provide explicit and computable formulas for their asymptotic variance that may be of interest to practitioners using these models. 
\end{itemize}
Besides the theoretical findings, we conduct numerical experiments to illustrate our results through concrete examples, including building confidence intervals in a simulation study and hypothesis testing on an Association of Tennis Professionals (ATP) dataset. 
\begin{table}[]
	\centering
	\begin{tabular}{cccc}
		\hline\hline
		Results &  Models & Comparison Graphs & Sparsity \\%& Elements \\
		\hline
		\cite{MR1724040} & BT & homogeneous &  $\Omega(1)$  \\
		\cite{MR2987494} & BT & heterogeneous & $\Omega(1)$   \\
		\cite{han2020asymptotic} & BT & homogeneous & $\widetilde{\Omega}(n^{-1/10})$   \\
		%\cite{MR4560195} & BT & homogeneous & ??? \\
		\cite{gao2023uncertainty} & BT & homogeneous & $\widetilde{\Omega}(n^{-1})$  \\
		This work & general  & heterogeneous & $\widetilde{\Omega}(n^{-1})$  \\
		\hline
	\end{tabular}\\
	
	\caption{Comparison of several asymptotic normality results in the pairwise comparison model literature. Homogeneous comparison graphs indicate the balanced data while heterogeneous graphs correspond to the imbalanced data.		
		The notation $\Omega(\cdot )$ is the standard Bachmann--Landau notation and $\widetilde{\Omega}(\cdot)$ means that the asymptotic relation holds up to polynomial terms of $\log n$.}
	\label{tab:my_label}
\end{table}

The rest of the paper is organized as follows.
In Section~\ref{sec:2}, we introduce the mathematical framework for general pairwise comparison models.
In Section~\ref{sec:3}, we establish the main results concerning asymptotic normality and non-asymptotic convergence rates (Theorems~\ref{CLT:1}, \ref{individual_error}, \ref{CI}) for the MLE and provide the relevant interpretations. 
In Section~\ref{sec:4}, we show that the MLE in many pairwise comparison models in the literature is asymptotically normal under minimal sparsity assumptions.
Furthermore, we offer explicit formulas for their asymptotic variance.
In Section~\ref{sec:5}, we validate the theoretical findings using both synthetic and real data.
In Section~\ref{flu}, we outline the key ideas for the proof of the main result, deferring the technical details to the appendix section.
In Section~\ref{sec:6}, we point out some future directions.

	\section{General pairwise comparison models}\label{sec:2}
We introduce the general pairwise comparison models proposed in \cite{han2022general}. 	
Consider an undirected graph $\cG = (\cV, \cE)$ with $n$ vertices, where $\cV = [n] := \{1, \ldots, n\}$ represents the subjects being modeled, and $\cE$ is the edge set of the comparison data, that is, for any $i, j \in \cV$, $(i, j) \in \cE$ if there exists a comparison between subjects $i$ and $j$. 
For ease of illustration, we assume that $\cG$ is simple so that there exists at most one edge between two subjects; the general case consisting of multiple edges can be considered similarly but requires heavier notation. 
We use $\delta\{i\} = \{j\in \cV: (i, j)\in \cE\}$ to denote the graph neighbourhood of $i$, which refers to the set of vertices adjacent to $i$.
For the moment, we assume that $\cG$ is known and describe the statistical model for the comparison outcomes. 

The comparison outcomes are modeled using latent scores. 
Let $\u^* = (u^*_1, \ldots, u^*_n)\in\R^n$ denote the latent score vector of the subjects under comparison. 
For instance, in a sports competition, one may think of $u^*_i$ as the internal strength of the $i$th player when competing with others. 
For $i, j\in [n]$, assuming there is a comparison between them, the outcome is a random variable $X_{ij}$ taking values in some symmetric subset $ \bbA\subseteq\R$.
The probability mass/density function of $X_{ij}$ follows a single-parameter family $f(x; u^*_i-u^*_j)$, depending on whether $X_{ij}$ is discrete or continuous.  
In the general pairwise comparison models, the function $f$ is assumed to be valid in the sense of the following definition.

\begin{definition}[Valid parameterization]\label{def:valid}
	A function $f: \bbA \times \mathbb{R} \to \mathbb{R}^+$, where $\bbA=-\bbA\subseteq\R$ denotes the possible comparison outcomes, is said to be valid if it satisfies the following assumptions:\\
	(A1).	 (normalization) For $y\in\R$, $\int_\bbA f(x;y) \d x = 1$, where the integral is interpreted as summation if $\bbA$ is discrete.\\
	(A2). (symmetry) $f$ is even with respect to $(x; y)$, that is, $f(x;y) = f(-x;-y) $, for $(x, y)\in \bbA\times\R.$\\
	(A3). (monotonicity) For $ x < 0 $, $ f(x;y) $ is decreasing in $y$, and $f(x;y)\to 0$ as $ y \to \infty $.\\
	(A4).  (boundedness) $\sup_{y\in\R}f(x; y)<+\infty$ for every $x\in \bbA$. \\
	(A5). (log-concavity) $f(x;y)$ is strictly log-concave with respect to $ y $.
\end{definition}
All pairwise comparison models mentioned in the introduction are special cases of general pairwise comparison models with suitable choices of valid $f$.
In the rest of the article, we shall always assume $f$ is valid.  

Having specified $f$, we can write down the likelihood function given the observed data and compute the MLE $\widehat{\u}$ for $\u^*$ under identifiability constraints; see Section \ref{sec:3}.  
It was shown in \cite{han2022general} that $\widehat{\u}$ is a uniformly consistent estimator for $\u^*$ under appropriate conditions on both $f$ and $\cG$. 
To identify graph configurations satisfying these conditions, it is often convenient to assume that $\cG$ is sampled from certain random graph ensembles. A commonly used one is the random graph model $\G\left(n, p_n, q_n\right)$ that generalizes the Erd\H{o}s--R\'enyi model \citep{MR0125031}. 
\begin{definition}[Random graph models] \label{Random_Graphs}
	$\G\left(n, p_n, q_n\right)$ is a random graph with vertices set $\cV=[n]$ where each edge $(i, j) \in \cV \times \cV, i \neq j$ is formed independently with deterministic probability $p_{i j, n} \in\left[p_n, q_n\right]$.
\end{definition}

\begin{remark}
	By definition, for fixed $p_n$ and $q_n$, different specifications of $\{p_{ij, n}\}_{1\leq i<j\leq n}$ may belong to the same $\G\left(n, p_n, q_n\right)$. 
	Such ambiguity is irrelevant to our discussion as only the upper and lower bounds are used in the analysis.  
\end{remark}

In the following section, we assume the comparison graph $\cG(\cV, \cE)$ is sampled from $\G\left(n, p_n, q_n\right)$. Meanwhile, for every $(i, j)\in \cE$, we assume that a single comparison outcome $X_{ij}\sim f(x; u^*_i-u^*_j)$ is observed. Note that $X_{ij} = -X_{ji}$, and we assume that $\{X_{ij}\}_{i<j, (i,j)\in\cE}$ are independent. It is worth noting that we allow $\lim\sup_{n \to \infty} q_n/p_n = +\infty$, which leads to heterogeneity of the observed data. In particular, one subject may have $nq_n$ comparisons while another has only $np_n$ comparisons. Such an imbalance of data volume will bring additional difficulty in theoretical data analysis.

\section{Main results}\label{sec:3}

This section establishes the asymptotic normality of the MLE $\widehat{\u}$ in the general pairwise comparison model introduced in Section \ref{sec:2}. 
In the following, we use $O(\cdot )$, $O_p(\cdot )$, $o(\cdot )$, and $o_p(\cdot)$ as a standard notation to represent the respective asymptotic order between two sequences. 

The log-likelihood function based on the observed outcomes $X_{ij}$ (conditional on $\cG(\cV, \cE)$) is 
\begin{align*}
	l(\u) = \sum_{(i, j)\in \cE}\log f(X_{ij}; u_i-u_j).
\end{align*}
Note that $l(\u)$ is invariant if all components of $\u$ are shifted by a constant.
To ensure model identifiability, we let $ {\bm 1}_n^\top\u = 0$, where $ {\bm 1}_n$ is the all-ones vector with the same dimension as $\u$. 
The MLE $\uu$ for $\u^*$ satisfies 
\begin{align}
	\uu = \argmax_{\u\in\R^{n}:  {\bm 1}_n^\top\u = 0}l(\u).\label{eq:mle}
\end{align}

To state the asymptotic normality result for $\widehat{\u}$, we need the following assumptions. 
Let $$M_n = \max_{i,j\in [n]}|u^*_i-u^*_j|$$
denote the dynamic range of $\u^*$.
The global discrepancy of the model, which measures the probability of a subject with the highest score winning over another with the lowest score, is defined as 
\begin{align*}
	c_{n,1} = \int_{\bbA\cap [0, +\infty)} f(x; M_n) \d x. 
\end{align*}
It can be verified using Definition~\ref{def:valid} that $c_{n,1}\in [1/2, 1)$. 

Let $g(x; y) = \partial_2(\log f(x; y)) = \partial_2 f(x; y)/f(x;y)$ be the Fisher score function, where $\partial_2$ denotes the partial derivative operator for the second argument $y$.  
The first assumption, which ensures the unique existence of the MLE, states that $c_{n,1}$ if converging to one, must have a controlled rate. 

\begin{assumption}\label{ass:0}
	The global discrepancy $c_{n,1}$ satisfies
	\begin{align}
		\frac{\log n}{np_n|\log c_{n,1}|}\to 0\qquad n\to\infty. \label{exist}
	\end{align}
\end{assumption}
%\vspace{0.2cm}

The next assumption is concerned with the tail decay of $g(X; y)$ with $X\sim f(x; y)$ for all $|y|\leq M_n$, where $g(X; y)$ is centered, that is, $\E[g(X; y)] = 0$. 

\begin{assumption}\label{ass:1}
	The sequence of random variables $\{g(X_y; y)\}_{|y|\leq M_n}$ is uniformly subgaussian, where $X_y\sim f(x; y)$. That is, 
	\begin{align*}
		c_{n, 2} = \max_{|y|\leq M_n}\|g(X_y; y)\|_{\psi_2}<\infty,
	\end{align*} 
	where $\|\cdot\|_{\psi_2}$ is the subgaussian norm \cite[Definition 2.5.6]{vershynin2018high}.
\end{assumption}
%\vspace{0.2cm}

In addition, we further requires a boundedness condition on the partial derivative of $g(x; y)$ in $y$. 

\begin{assumption}\label{ass:2}
	The partial derivative of the score function $g$ for $y$ is uniformly bounded from both above and below:
	\begin{align*}
		0<c_{n,3} = \inf_{x\in \bbA, |y|\leq M_n+1}|\partial_2 g(x; y)|\leq \sup_{x\in \bbA, |y|\leq M_n+1}|\partial_2 g(x; y)| = c_{n,4}<\infty.
	\end{align*}
\end{assumption}

Assumptions \ref{ass:0}-\ref{ass:2} are sufficient to guarantee the uniform consistency of the maximum likelihood estimator \citep{han2022general}.
The common pathway from consistency to asymptotic normality is via Taylor's expansion. As a result, one often require more stringent regularity conditions on higher-order derivatives.
%In our case, we need to work with the (pseudo)-inverse of the Hessian matrix to build the connection with the central limit theorem.  %{\color{red}The following assumption is a mild sufficient condition ensuring that the second smallest eigenvalue of the Hessian matrix is not close to zero.}

\begin{assumption}\label{ass:3}
	The second-order partial derivative of the score function $g(x, y)$ for $y$ is uniformly bounded for all $|y|\leq M_n+1$:
	\begin{align*}
		c_{n,5}= \sup_{x\in \bbA, |y|\leq M_n+1}|\partial_{22}g(x; y)|<\infty.
	\end{align*}
\end{assumption}

\begin{theorem}[Asymptotic normality]\label{CLT:1}
	Let $\cG\sim\G(n, p_n, q_n)$.
	Under Assumptions \ref{ass:0}--\ref{ass:3}, 
	if
	\begin{align}
		\beta_n = \max\left\{ \frac{c_{n,2}^2 c_{n,4}^{5/2}c_{n,5}}{c_{n,3}^{5}}, \frac{c_{n,2} c_{n,4}^{11/2}}{c_{n,3}^{6}}\right\}\left\{\frac{q_n^{10}(\log n)^8}{np_n^{11}}\right\}^{1/2}\to 0\qquad n\to\infty, \label{good}
	\end{align}
	then for each $i\in [n]$, $\{\rho_i(\u^*)\}^{-1/2}(\widehat{u}_i - u_i^*) \to \mathcal{N}(0, 1)$, where $\mathcal{N}(0,1)$ is the standard normal distribution and the asymptotic variance $\rho_i(\u^*)$ is given by
	\begin{align}
		\rho_i(\u^*) 
		&= \left[\sum_{j\in\delta\{i\}}\int_{\bbA} \frac{\{\partial_2 f(x; u^*_i-u^*_j)\}^2}{f(x; u^*_i-u^*_j)}\d x\right]^{-1}.\label{avc}
	\end{align}
	Moreover, for any finite $S = \{i_1, \ldots, i_s\}\subset \N$ where $S$ is independent of $n$, $\{\rho_{i_1}(\u^*)\}^{-1/2}(\widehat{u}_{i_1} - u_{i_1}^*), \ldots, \{\rho_{i_s}(\u^*)\}^{-1/2}(\widehat{u}_{i_s} - u_{i_s}^*)$ are asymptotically independent.  
\end{theorem}
\begin{remark}
	Subject to additional conditions, such as a uniform upper bound on the number of edges between two subjects that is independent of $\u^*$, Theorem \ref{CLT:1} still holds when two subjects are compared multiple times.   
	Under such circumstances, the asymptotic variance formula \eqref{avc} remains valid if one counts multiple edges for each adjacent vertex. 
\end{remark}

A proof sketch for Theorem \ref{CLT:1} is provided in Section \ref{flu}. To better understand Theorem~\ref{CLT:1}, note that although the multiplicative factor involving $c_{n, i}$ may seem daunting, it is $O(1)$ under additional regularity and boundedness conditions.

\begin{corollary}\label{myco}
	Assume that $M^*=\sup_{n}M_n<\infty$ and $\bbA$ is finite.
	If $f(x;y)$ is a vaild function, and $g(x;y), \partial_2 g(x;y), \partial_{22}g(x;y)$ are continuous in $y$, then the asymptotic normality result in Theorem  \ref{CLT:1} holds if $q^{10}_n(\log n)^8/(np^{11}_n)\to 0$ as $n\to\infty$.  
\end{corollary}

Corollary~\ref{myco} can be further simplified under balanced data conditions. If $\sup_{n}q_n/p_n<\infty$, then the condition in Corollary~\ref{myco} reduces to $(\log n)^8/(np_n)\to 0$ as $n\to\infty$. 
This bound coincides with the well-known Erd\H{o}s--R\'enyi connectivity threshold up to logarithmic factors and is thus close to optimal. In general, Corollary~\ref{myco} allows $p_n$ and $q_n$ to have different asymptotic orders. For instance, if $p_n = n^{-b}$ for $0\leq b<1$, then choosing $q_n = n^{-b'}$ for any $b'\geq 0$ with $b-(1-b)/10<b'\leq b$ satisfies the required conditions. To the best of our knowledge, this is the first result that works for general comparison models with imbalanced data.

Theorem \ref{CLT:1} provides the limiting distribution of the MLE, which lays the foundation for further statistical inference, such as the construction of interval-based estimates, testing hypothesis on the difference of merits, and testing the hypothesis on ranking position \citep{doi:10.1080/01621459.2024.2316364}. Additionally, we can derive a non-asymptotic error rate for each parameter as well as a practical construction of interval-valued estimators, which are summarized in Theorem \ref{individual_error} and Theorem \ref{CI}, respectively. 

\begin{theorem}[Individual estimation error]\label{individual_error}
	Let $\cG\sim\G(n, p_n, q_n)$. If both Assumptions \ref{ass:0}--\ref{ass:3} and \eqref{good} hold, then there exists an absolute constant $C>0$ such that for each $i\in [n]$, the following holds for all sufficiently large $n$ with probability tending to one:
	\begin{equation}\label{individual}
		|\widehat{u}_{i} - u_{i}^*| \leq C\cdot\frac{c_{n,2}}{c_{n,3}}\left(\frac{\log n}{|\delta\{i\}|}\right)^{1/2}.
	\end{equation}
\end{theorem}
We call the result in \eqref{individual} the individual estimation error bound. Compared with the existing works \citep{han2020asymptotic,chen2022optimal}, this result presents new findings on the individual estimation error. Specifically, previous works considered the overall estimation error $\max_{i \in[n]} |\widehat{u}_{i} - u_{i}^*|$, while our result in \eqref{individual} further illustrates the error for each individual. As expected, this error depends on the degree $|\delta\{i\}|$, which serves as a local normalization factor for each individual. When the comparison graph is homogeneous, individual estimation errors converge at the same rate, matching the results in \cite{han2020asymptotic} and \cite{chen2022optimal}. For heterogeneous graphs, however, individuals will have different convergence rates, which cannot be reflected in the worst-case error $\max_{i \in[n]} |\widehat{u}_{i} - u_{i}^*|$. 

Next, we construct the confidence interval estimators. Since the asymptotic variance $\rho_i(\u^*) $ is unknown, it is natural to consider the following plug-in variance estimator:
\begin{align}
	\rho_i(\uu) 
	&= \left[\sum_{j\in\delta\{i\}}\int_{\bbA} \frac{\{\partial_2 f(x; \widehat{u}_{i}-\widehat{u}_{i})\}^2}{f(x; \widehat{u}_{i}-\widehat{u}_{i})}\d x\right]^{-1}.\label{estimate_var}
\end{align}
Then for any fixed $i \in [n]$ and $\alpha \in (0,1)$, the $(1-\alpha)$-confidence interval  of $u_i^*$ is
\begin{align*}
	\mathcal{CI}_i = 	\big(\widehat{u}_{i}-z_{\alpha/2}\{\rho_{i}(\uu)\}^{-1/2}, \widehat{u}_{i} + z_{\alpha/2}\{\rho_{i}(\uu)\}^{-1/2} \big),
\end{align*}
where  $z_{\alpha/2}$ is the  $(1- \alpha/2)$-quantile of standard normal distribution. 
Note that confidence intervals constructed using consistently estimated plug-in estimators are often automatically valid due to Slutsky's theorem if the parameters to be estimated are independent of $n$. In our case, the increasing dimension of the parameter vector makes such a result less obvious, though it remains true, as justified by the following theorem.
\begin{theorem}[Confidence interval]\label{CI}
	Let $\cG\sim\G(n, p_n, q_n)$. Under Assumptions \ref{ass:0}--\ref{ass:3}, if \eqref{good} holds, then  for any fixed $i \in [n]$ and $\alpha \in (0,1)$, 
	\begin{align*}
		\mathbb{P}\left(u^*_i \in \mathcal{CI}_i   \right) \to 1-\alpha  \ \text{ as } n \to \infty.
	\end{align*}
\end{theorem}
The computation of $\mathcal{CI}_i$ is straightforward after obtaining $\uu$. 
For example, if one wants to construct a $95\%$-confidence interval for the $i$th subject, one can take $\mathcal{CI}_i = \big(\widehat{u}_{i}-1.96\{\rho_{i}(\uu)\}^{-1/2}, \widehat{u}_{i} + 1.96\{\rho_{i}(\uu)\}^{-1/2} \big)$.

\section{Examples}\label{sec:4}

We establish the asymptotic normality of the MLEs for a variety of pairwise comparison models found in the existing literature under conditions characterized by $p_n, q_n$, and $M_n$. These models encompass the BT model, the Thurstone--Mosteller model, the Rao--Kupper model, the Davidson model, the cumulative link model (with four outcomes), and the paired cardinal model. In particular, explicit formulas are derived for their asymptotic variance.
For conveniences, we let $\Delta_{ij} = u^*_i - u^*_j$ for $i, j\in [n]$. 

\begin{example}[BT model]
	The BT model assumes a binary comparison outcome between two subjects and parameterizes using a logistic link function:
	\begin{align*}
		f(1; y) = \frac{e^y}{1+ e^y}; \quad f(-1 ; y) = \frac{1}{1+ e^y}.
	\end{align*} 
	It can be verified that both conditions in \eqref{exist} and \eqref{good} are satisfied if
	\begin{align*}
		e^{6M_n}\left\{\frac{q_n^{10}(\log n)^8}{np_n^{11}}\right\}^{1/2}\to 0\qquad n\to\infty.
	\end{align*}	
	The asymptotic variance $\rho_i(\u^*)$ for the subject $i \in [n]$ is given by
	\begin{align*}
		\rho_i(\u^*) = \left\{\sum_{j \in \delta\{i\}} \frac{ e^{\Delta_{ij}} }{\left(1+ e^{\Delta_{ij}} \right)^2}\right\}^{-1}.
	\end{align*}
	That is, $\{{\rho_i(\u^*)}\}^{-1/2} (\widehat{u}_i - u^*_i)\to \mathcal{N}(0,1)$ as $n\to\infty$.\\
\end{example}

\begin{example}[Thurstone--Mosteller model]
	The Thurstone--Mosteller model considers an alternative link function than the BT model:
	\begin{align*}
		f(1; y) = \Phi(y); \quad f(-1 ; y) = 1-\Phi(y),
	\end{align*} 
	where $\Phi(y) = \int_{-\infty}^y \varphi(z) \d z$ with $\varphi(z) = (2\pi)^{-1/2}e^{-z^2/2}$.
	It can be verified that both conditions in \eqref{exist} and \eqref{good} are satisfied if
	\begin{align*}
		M_n^{6}e^{3M^2_n}\left\{\frac{q_n^{10}(\log n)^8}{np_n^{11}}\right\}^{1/2}\to 0\qquad n\to\infty.
	\end{align*}			
	The asymptotic variance $\rho_i(\u^*)$ for the subject $i \in [n]$ is given by
	\begin{align*}
		\rho_i(\u^*) = \left \{\sum_{j \in \delta\{i\}} \frac{\varphi^2(\Delta_{ij})}{\Phi(\Delta_{ij})(1-\Phi(\Delta_{ij}))}\right \}^{-1}.
	\end{align*}
	That is, $\{{\rho_i(\u^*)}\}^{-1/2} (\widehat{u}_i - u^*_i)\to \mathcal{N}(0,1)$ as $n\to\infty$.\\
\end{example}

\begin{remark}
Both the BT and Thurstone--Mosteller model models can be viewed as special cases of the Random Utility Model framework, which was formally developed for discrete choice analysis by \cite{mcfadden1973conditional}. Specifically, the BT model arises from the assumption of i.i.d. Gumbel-distributed random utilities, while the Thurstone--Mosteller model assumes they are normally distributed.  Our central findings can be generalized to the other types random utility models, which satisfies (A1)--(A5) in Definition \ref{def:valid}.
\end{remark}

\begin{example}[Rao--Kupper model]
	The Rao--Kupper model extends the BT model by incorporating ties. 
	The link function is given by
	\begin{align*}
		f(1 ; y)&=\frac{e^y}{e^y+\theta}; \quad f(0 ; y)=\frac{\left(\theta^2-1\right) e^y}{\left(e^y+\theta\right)\left(\theta e^y+1\right)}; \quad		f(-1 ; y)=\frac{1}{\theta e^y+1},
	\end{align*}
	where $\theta>1$ is the threshold parameter which is predetermined. 	
	It can be verified that both conditions in \eqref{exist} and \eqref{good} are satisfied if
	\begin{align*}
		e^{6M_n}\left\{\frac{q_n^{10}(\log n)^8}{np_n^{11}}\right\}^{1/2} \to 0\qquad n\to\infty.
	\end{align*}		
	The asymptotic variance $\rho_i(\u^*)$ for the subject $i \in [n]$ can be computed as
	\begin{align*}
		\rho_i(\u^*)=	\left\{\sum_{j \in \delta\{i\}}\left[\frac{\theta^2 e^{\Delta_{ij}}}{(\theta+ e^{\Delta_{ij}} )^3}+\frac{\theta^2(\theta^2-1)e^{\Delta_{ij}}(1-e^{2\Delta_{ij}})^2}{( e^{\Delta_{ij}} +\theta)^3(\theta e^{\Delta_{ij}} +1)^3}+\frac{\theta^2e^{2\Delta_{ij}}}{(\theta e^{\Delta_{ij}} +1)^3}\right]\right\}^{-1}.
	\end{align*}
	That is, $\{{\rho_i(\u^*)}\}^{-1/2} (\widehat{u}_i - u^*_i)\to \mathcal{N}(0,1)$ as $n\to\infty$.\\
\end{example}

\begin{example}[Davidson model]
	As opposed to the Rao--Kupper model, the Davidson model considers an alternative parameterization of being tied:
	\begin{align*}
		f(1 ; y)=\frac{e^y}{e^y+\theta e^{\frac{y}{2}}+1}; \quad f(0 ; y)=\frac{\theta e^{\frac{y}{2}}}{e^y+\theta e^{\frac{y}{2}}+1}; \quad
		f(-1 ; y)&=\frac{1}{e^y+\theta e^{\frac{y}{2}}+1},
	\end{align*}
	where $\theta>0$ is assumed to be prefixed.
	It can be verified that both conditions in \eqref{exist} and \eqref{good} are satisfied if
	\begin{align*}
		e^{{3M_n}}\left\{\frac{q_n^{10}(\log n)^8}{np_n^{11}}\right\}^{1/2}\to 0\qquad n\to\infty.
	\end{align*}		
	The asymptotic variance $\rho_i(\u^*)$ for the subject $i \in [n]$ is given by
	\begin{align*}
		\rho_i(\u^*) = \left\{\sum_{j \in \delta\{i\}}\frac{e^{\Delta_{ij}}(\theta e^{\Delta_{ij}/2}+2)^2+\theta e^{\Delta_{ij}/2}(1-e^{\Delta_{ij}})^2 +(2e^{\Delta_{ij}}+\theta e^{\Delta_{ij}/2})^2}{4( e^{\Delta_{ij}}  + \theta e^{\Delta_{ij}/2} +1)^3}\right\}^{-1}.
	\end{align*}
	That is, $\{{\rho_i(\u^*)}\}^{-1/2} (\widehat{u}_i - u^*_i)\to \mathcal{N}(0,1)$ as $n\to\infty$. \\
\end{example}

\begin{example}[Cumulative link model with four outcomes] \label{Example_5}
	The cumulative link model is defined via the ordinalization of some distribution function $F$ \citep{agresti1992analysis}. 
	The BT model, the Rao--Kupper model, and the Thurstone--Mosteller model can all be considered as specific cases of the cumulative link model with appropriate choices of $F$. 
	In this example, we focus on the situation when $F$ is logistic with four outcomes, which will be used later in Section \ref{sec:5}.    
	In this case, $\mathbb A=\{-2,-1,1,2\}$ and the corresponding link function takes the form of 
	\begin{align}
		f(1;y)=\frac{(\theta-1) e^{y}}{\left(\theta+e^{y}\right)\left(1+e^{y}\right)}, \ \ \ f(2;y)=\frac{e^{y}}{\theta+e^{y}},\label{tah}
	\end{align} 
	where $\theta>1$ is the threshold parameter that is predetermined. It can be verified that both conditions in \eqref{exist} and \eqref{good} are satisfied if
	\begin{align*}
		e^{6M_n}\left\{\frac{q_n^{10}(\log n)^8}{np_n^{11}}\right\}^{1/2}\to 0\qquad n\to\infty.
	\end{align*}		
	The asymptotic variance $\rho_i(\u^*)$ for the subject $i \in [n]$ is given by
	{\small \begin{align*}
			{\rho_i(\u^*) =	\left\{\sum_{j \in \delta\{i\}}\left[\frac{\theta^2 e^{\Delta_{ij}}}{(\theta+ e^{\Delta_{ij}} )^3}+\frac{\theta^2 e^{2\Delta_{ij}}}{(\theta e^{\Delta_{ij}}+1 )^3}+\frac{(\theta-1)e^{\Delta_{ij}}}{(1+e^{\Delta_{ij}})^3}\left(\frac{(\theta-e^{\Delta_{ij}})^2}{(\theta+e^{\Delta_{ij}})^3}+\frac{1-\theta e^{2\Delta_{ij}}}{(\theta e^{\Delta_{ij}}+1)^3}\right)\right]\right\}^{-1}.}
	\end{align*}}
	\noindent That is, $\{{\rho_i(\u^*)}\}^{-1/2} (\widehat{u}_i - u^*_i)\to \mathcal{N}(0,1)$ as $n\to\infty$. 
	The result above can be extended to the cumulative link model with the logistic link function and finite outcomes but involve more complicated formulas for asymptotic variances. 
\end{example}

\begin{example}[Paired cardinal model]
	The paired cardinal model can be viewed as a fully observed version of the Thurstone--Mosteller model, where the comparison outcome takes values in $\mathbb A=\R$, and the corresponding link function is parameterized using the normal density:
	\begin{align*}
		f(x; y) = \left(2\pi\sigma^2\right)^{-\frac12}e^{-\frac{(x-y)^2}{2\sigma^2}},
	\end{align*}
	where $\sigma>0$ is assumed to be prefixed.
	It can be verified that both conditions in \eqref{exist} and \eqref{good} are satisfied if
	\begin{align*}
		\max\left\{M_ne^{\frac{M_n^2}{2\sigma^2}}\frac{\log n}{np_n}, \left({\frac{q_n^{10}(\log n)^8}{np_n^{11}}}\right)^{1/2}\right\}\to 0\qquad n\to\infty.
	\end{align*}
	The asymptotic variance $\rho_i(\u^*)$ for the subject $i \in [n]$ is given by
	\begin{align*}
		\rho_i(\u^*) = \frac{\sigma^2}{|\delta\{i\}|},
	\end{align*}
	which is independent of $\u^*$. 
	That is, $\{{\rho_i(\u^*)}\}^{-1/2} (\widehat{u}_i - u^*_i)\to \mathcal{N}(0,1)$ as $n\to\infty$. 
\end{example}

	\section{Numerical studies}\label{sec:5}

\subsection{Synthetic data}
We conduct a simulation study to verify the asymptotic normality of the MLE. We focus on three comparison models, the Davidson model, the Rao--Kupper model, and the paired cardinal model, as discussed in the preceding section.

\begin{figure}[htb]
	\centering
	\includegraphics[width=0.328\textwidth, height=1.7in, trim={1.2cm 0.5cm 1.5cm 1cm},clip]{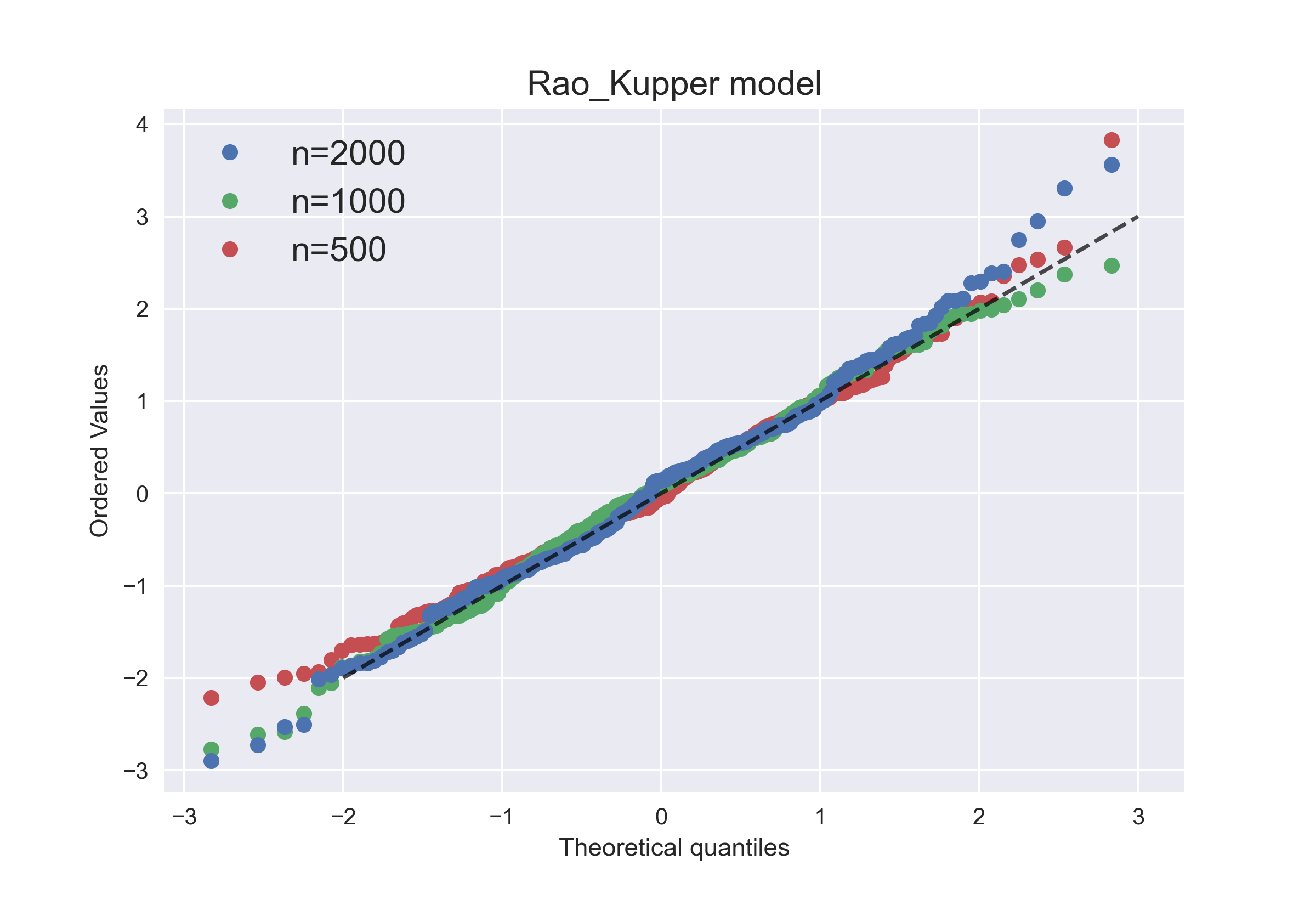}
	\includegraphics[width=0.328\textwidth, height=1.7in, trim={1.2cm 0.5cm 1.5cm 1cm},clip]{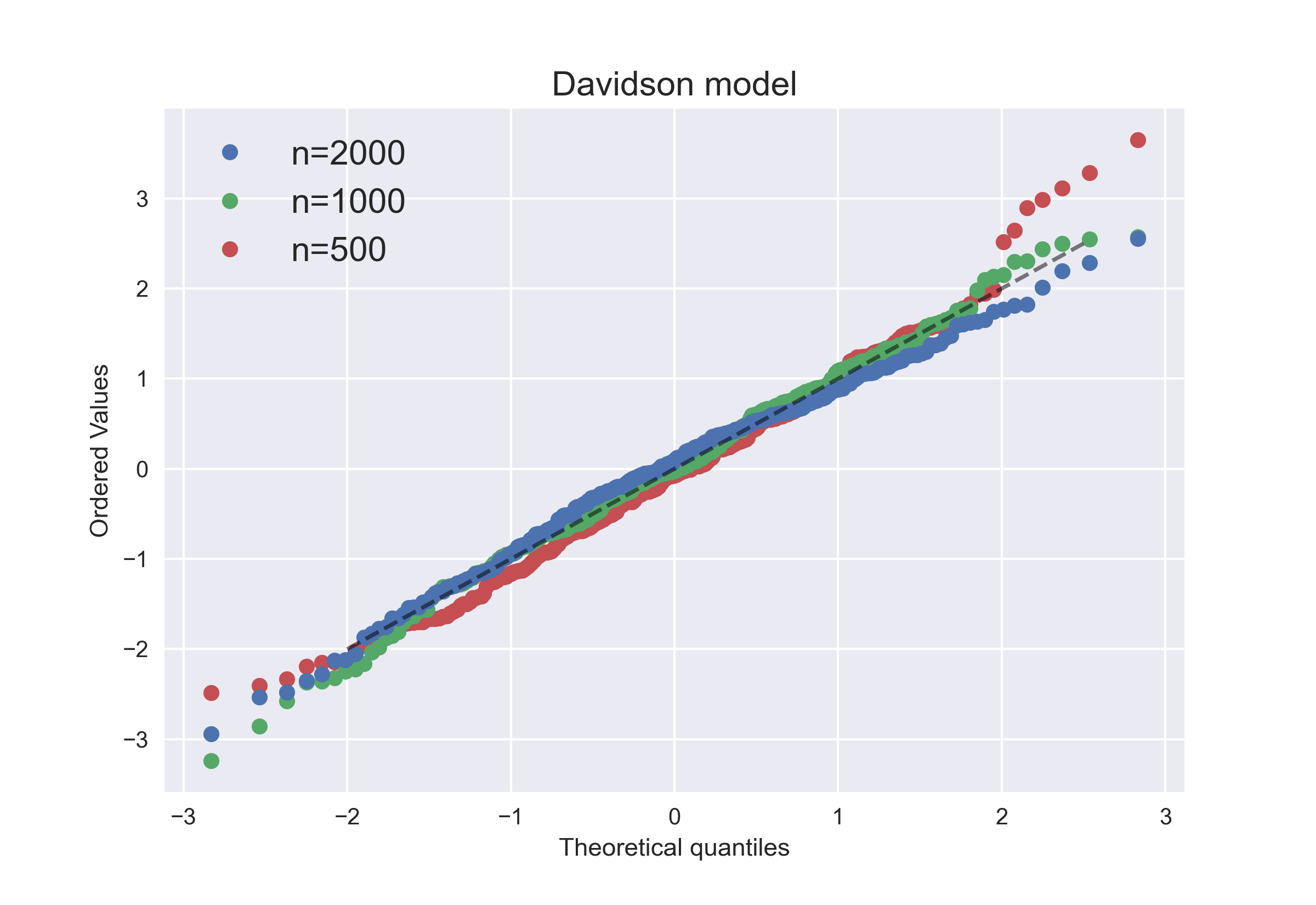}
	\includegraphics[width=0.328\textwidth, height=1.7in, trim={1.2cm 0.5cm 1.5cm 1cm},clip]{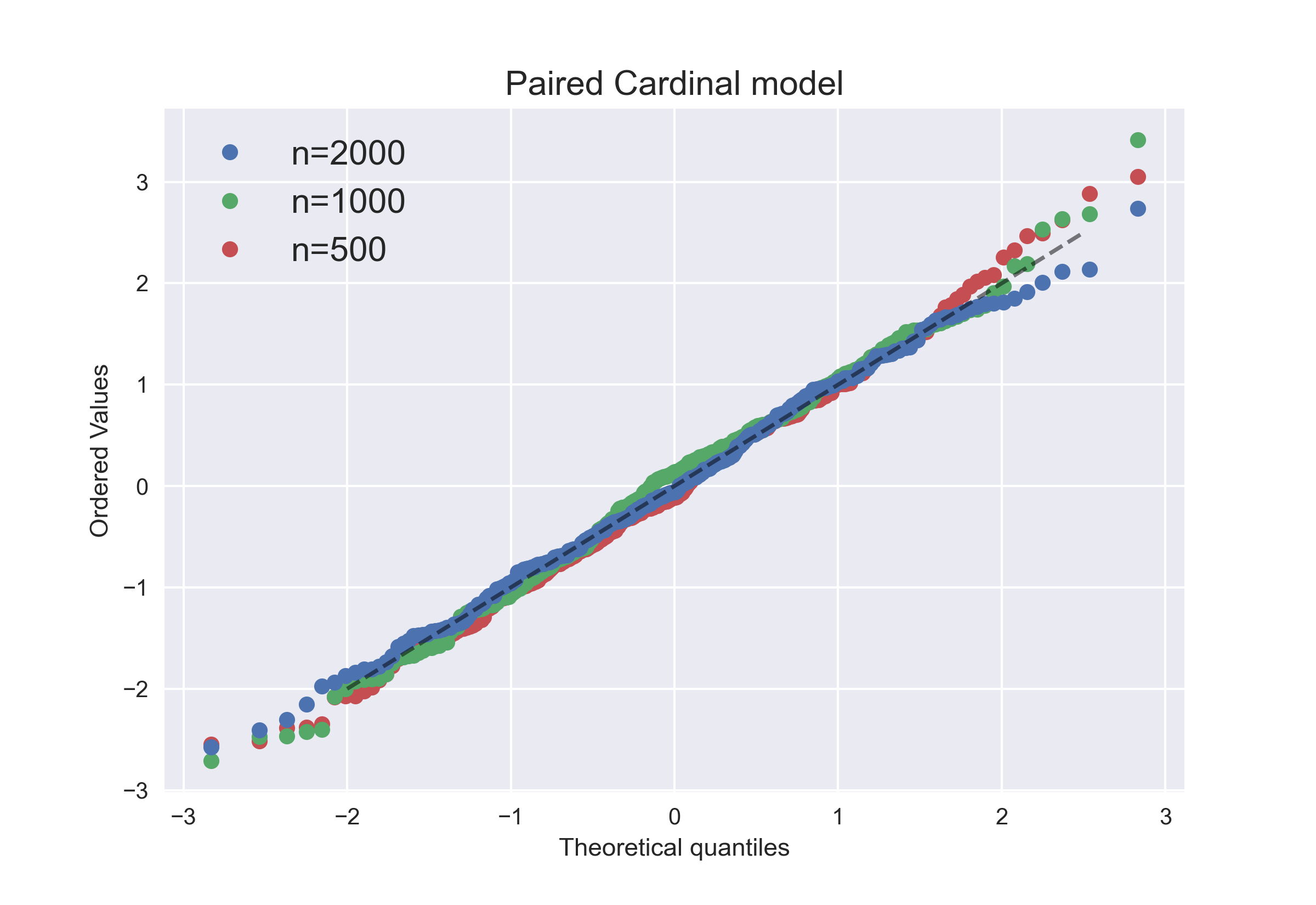}		
	\caption{Quantile-quantile plots comparing the theoretical and sample
		quantiles of $\widehat{u}_1$ in the Rao--Kupper model, the Davidson model, and the paired cardinal model with fixed $M_n=1$. Simulation parameters: For each $n$, the empirical
		quantile curves are based on 300 experiments. }
	\label{Fig.main.1}
\end{figure}

We utilize the general random graph model introduced in Definition \ref{Random_Graphs} to generate comparison graphs, where we set $p_n=n^{-1/2}$ and $q_n=p_n (\log n) $, with the individual edge probabilities chosen uniformly at random between them. The total number of subjects $n$ in all the tested models is chosen from the set $\{500, 1000, 2000\}$. To generate the utility vector $\u^*$, we uniformly select values from the range $[-M_n, M_n]$, where the dynamic range $M_n$ is chosen from $\{1, \log(\log n)\}$. For global parameters in each model, we set $\theta=1$ in the Davidson model, $\theta=2$ in the Rao--Kupper model, and $\sigma=2$ in the paired cardinal model. Each simulation scenario is characterized by the pair $(n, M_n)$ and repeated $300$ times. We report the average standard deviation of the estimated parameters and the coverage probability of the $95\%$ confidence interval, as shown in Table \ref{table1}. We also examine the $z$-scores of the first coordinate of ${\uu}$ over 300 experiments and plot the corresponding quantile-quantile plots against normal quantiles in Figure \ref{Fig.main.1}.

In Table \ref{table1}, we observe a decreasing trend in the average standard deviation as $n$ increases, which is consistent with our theoretical findings on asymptotic variance in \eqref{avc}. Additionally, we note variations in standard deviation relative to $M_n$ across different models. In particular, the standard deviations are increasing in $M_n$ in the Rao--Kupper and Davidson models, while remaining unchanged in the paired cardinal model. This observation aligns with their respective asymptotic variance formulas $\rho_i(\u^*)$ computed in Section \ref{sec:4}; in particular, $\rho_i(\u^*)$ in the paired cardinal model is \emph{independent} of $M_n$.

The empirical coverage probabilities of the constructed confidence intervals for all models under comparison closely approximate the desired coverage of $95\%$. This, combined with the diagonal alignment of the data in Figure \ref{Fig.main.1}, suggests that the expected asymptotic normality holds. 
\begin{table}[htb]
	\caption{	\centering Summarized simulation results over 300 replications.}
	\label{table1}
	\centering
	\begin{tabular}{ccccc|ccc}
		\hline 
		& & \multicolumn{3}{c|}{$M_n = 1$} & \multicolumn{3}{c}{$M_n = \log (\log n)$}  \\ \hline
		\multicolumn{2}{c}{$n$} & 500 & 1000 & 2000 & 500 & 1000 & 2000 \\ \hline
		\multicolumn{1}{c}{Rao--Kupper}& Standard deviation & $0.157$ & $0.126$ & $0.101$ & $0.172$ & $0.139$ & $0.113$ \\ %\hline
		& Coverage probability & $0.949$ & $0.949$ & $0.949$ & $0.948$ & $0.949$ & $0.949$ \\ \hline
		\multicolumn{1}{c}{Davidson}& Standard deviation &$0.209$ & $0.167$ & $0.134$ & $0.225$ & $0.181$ & $0.147$ \\ %\hline
		& Coverage probability & $0.949$ & $0.950$ & $0.950$ & $0.949$ & $0.951$ & $0.950$ \\ \hline
		\multicolumn{1}{c}{Paired cardinal}& Standard deviation &$0.164$ & $0.131$ & $0.105$ & $0.163$ & $0.131$ & $0.105$ \\ %\hline
		& Coverage probability & $0.950$ & $0.950$ & $0.950$ & $0.950$ & $0.950$ & $0.949$ \\ \hline %\hline
	\end{tabular}
\end{table}

\subsection{Real data example}
\begin{figure}[htp]
\begin{subfigure}{0.45\textwidth}
	\centering
	\fbox{\includegraphics[width=\textwidth, height=2in,clip]{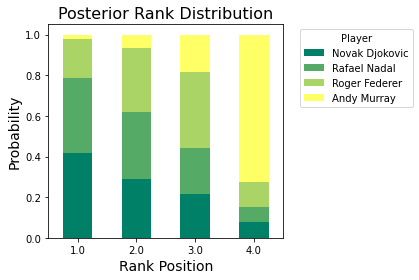}}
	\caption{Posterior ranking probability histogram.}
	\label{Fig.sub.21}
\end{subfigure}
	\begin{subfigure}{0.45\textwidth}
		\centering
		\fbox{\includegraphics[width=\textwidth, height=2in,clip]{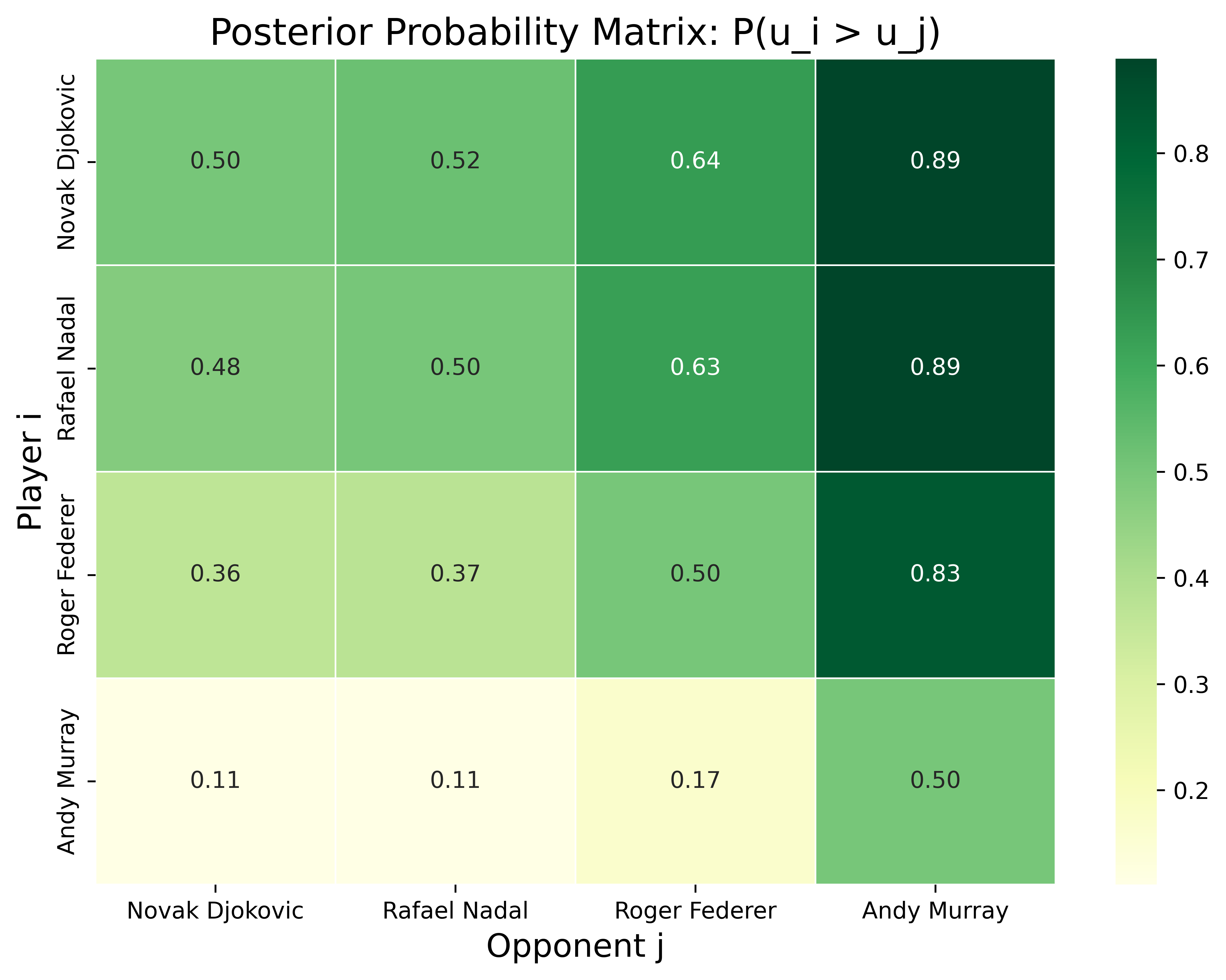}}
		\caption{Win probability heatmap.}
		\label{Fig.sub.22}
	\end{subfigure}
	\caption{ The posterior ranking and pairwise win probability for each player.} 
	\label{Fig.main.2}
\end{figure}

\begin{figure}[htp]
\begin{subfigure}{0.45\textwidth}
	\centering
		\fbox{\includegraphics[width=\textwidth, height=2in,clip]{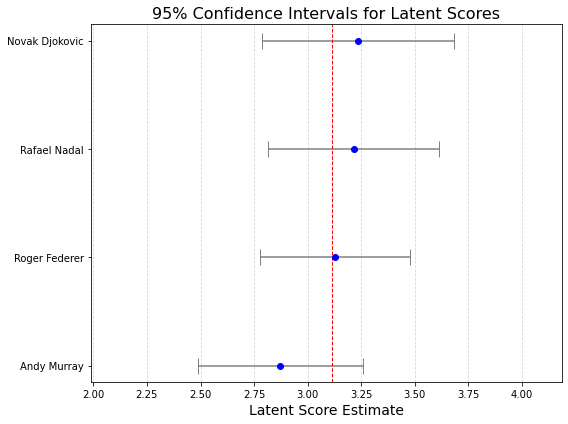}}
	\caption{The $95\%$ confidence interval of latent score.}
	\label{Fig.sub.31}
\end{subfigure}
\begin{subfigure}{0.45\textwidth}
	\centering
\fbox{\includegraphics[width=\textwidth, height=2in,clip]{pics/compete_CI.png}}
	\caption{The $95\%$ confidence interval of difference.}
	\label{Fig.sub.32}
\end{subfigure}
	\caption{ The $95\%$ confidence interval of  the latent score of individual players and the pairwise difference between players' scores.} 
	\label{Fig.main.3}
\end{figure}

{
In the realm of men's tennis, the dominance of Novak Djokovic, Roger Federer, Rafael Nadal, and Andy Murray, collectively known as the Big Four, has undeniably dominated the first quarter of the 21st century. However, recent discussions have emerged, advocating for the exclusion of Andy Murray in what could potentially become the Big Three. 
In this section, we consider a statistical approach to model the ATP data and uncover the latent scores that provide evidence in favour of the Big Four.
We consider modeling the comparison data from ATP\footnote{www.tennis-data.co.uk} 
using the cumulative link model with four outcomes from Example \ref{Example_5}. The ATP dataset collects the comparison results of each game from 1968 to 2016 by 2742 players. In this study, we concentrate on the best-of-3 matches, which result in four possible outcomes $\{\text{0:2, 1:2, 2:1, 2:0} \}$, represented as $\bbA=\{-2,-1,1,2\}.$
In particular, transferring the ordinal outcomes \eqref{tah} into the current setting reads
$$
\begin{aligned}
	&f(X_{ij}=1,u_i-u_j)=\frac{(\theta-1) e^{u_i-u_j}}{\left(\theta+e^{u_i-u_j}\right)\left(1+e^{u_i-u_j}\right)}, \\
	& f(X_{ij}=2,u_i-u_j)=\frac{e^{u_i-u_j}}{\theta+e^{u_i-u_j}}, \ \ \ \theta>1 .
\end{aligned}
$$

 The estimated latent scores $\widehat{u}_i$ and corresponding standard deviation $(\widehat{\rho}_i)^{1/2} $ of the 4 candidates are presented in Table \ref{ATP}, and it shows that Djokovic, Nadal, and Federer exhibit comparable performance levels, while Murray's score is slightly lower. Moreover, the overlap in standard deviations implies potential ranking uncertainty. 
 
To  quantify the ranking uncertainty, %we generate a posterior rank distribution based on  $\widehat{u}_i$ and $(\widehat{\rho}_i)^{1/2} $ using parametric bootstrap. Specifically, 
we use parametric bootstrap to simulate 10,000 samples by the asymptotic normal distributions $\mathcal{N}(\widehat{u}_i, \widehat{\rho}_i)$. Then we calculate empirical ranks across Monte Carlo samples and construct histogram for each rank position, shown in Figure \ref{Fig.sub.21}, where each bar indicates the probability of each player ranking from 1st through 4th. Figure \ref{Fig.sub.21} reveals that Djokovic has the highest probability of achieving the top rank, closely followed by Nadal and Federer, while Murray consistently ranks last, indicating the gap between him and the Big Three.
 
 To further investigate pairwise competitiveness among players, we construct a posterior probability matrix. It summarizes the likelihood that player $i$ outperforms player $j$ based on their estimates. Specifically, we compute:
 $$
 \P\left(u_i>u_j\right)=\mathbb{P}\left(Z_i-Z_j>0\right),
 $$
 where $Z_i \sim \mathcal{N}\left(\widehat{u}_i, \widehat{\rho}_i\right)$ and $Z_j \sim \mathcal{N}\left(\widehat{u}_j, \widehat{\rho}_j\right)$ are independently drawn from the posterior distributions of players $i$ and $j$, respectively.
 The asymptotic distribution guarantees that
 $$
 \mathbb{P}\left(u_i>u_j\right)=\Phi\left(\frac{\widehat{u}_i-\widehat{u}_j}{\sqrt{\widehat{\rho}_i+\widehat{\rho}_j}}\right),
 $$
 where $\Phi(\cdot)$ denotes the standard normal cumulative distribution function. The heatmap in Figure \ref{Fig.sub.22} displays the resulting posterior win probability matrix , where each entry 
 $(i,j)$ shows the empirical probability that player $i$ wins player $j$ in terms of latent score under 10,000 repetitions. We can observe that Djokovic and Nadal have approximately an equal chance of outperforming each other, and both have $60$-$65\%$ probability of outperforming Federer. Murray has the lowest win probabilities against all other players, consistent with his lower latent score. In general,  Djokovic and Nadal are the top players with extremely close latent strength, Federer is a strong challenger but slightly behind and Murray most frequently ranked the lowest.
 
 To further assess the uncertainty in individual scores, we construct confidence intervals for latent scores.  Figure \ref{Fig.sub.31} presents the $95\%$ confidence intervals for $\widehat{u}_i$ of each player under the CLM4 model. We use standard normal approximation to the asymptotic distribution to construct the confidence intervals, i.e., $\widehat{u}_i \pm z_{0.95} \sqrt{\widehat{\rho}_i}$. Figure \ref{Fig.sub.31} confirms the dominance of Djokovic and Nadal, the competitive strength of Federer, and a notable performance gap for Murray. However, the overlapping of all confidence intervals shows that these 4 players are not significantly different. To further validate the distinction with statistically significance, we construct the $95\%$ confidence interval for latent score difference. To construct the intervals, we apply the Bonferroni-adjusted critical value to the pairwise normal differences
 $$
 \widehat{u}_i-\widehat{u}_j \pm z_{1-\alpha /(2 K)} \cdot \sqrt{\widehat{\rho}_i+\widehat{\rho}_j}
 $$
 where $K$ is the number of pairwise comparisons and $z_{1-\alpha /(2 K)}$ is the Bonferroni-corrected quantile. Figure \ref{Fig.sub.32} shows that  none of the intervals exclude point zero, indicating that although Murray seems rank the least from previous experiments, all four players are not significantly different from others at the $95\%$ confidence level.
  
 In summary, both point estimates and posterior analyses indicate a discernible ranking between Murray and the Big Three. However, the simultaneous confidence intervals reveal that none of the pairwise differences are statistically significant. This highlights the critical role of formal statistical inference in practical applications.}

%%%%%%%%%%%%%%%%%%%%%%%%%%%%%%%%%%%%%%%%%%%%%%%%%%%%%%%%%%%%%%%%%%%%

%%% comment the following to accelerate compile
\begin{table}[!h]
	\caption{The estimated scores and their standard deviation of four players by using CLM4. }
	\label{ATP}
	\centering
	\begin{tabular}{ccc}
		Players & Estimated scores $\widehat{u}_i$ & Estimated standard deviation $(\widehat{\rho}_i)^{1/2} $  \\ \hline
		Novak Djokovic	&3.235&	$0.229$\\
		Rafael Nadal&	3.214&	$0.203$\\
		Roger Federer&	3.129&	$0.179$\\
		Andy Murray&	2.872	&$0.196$ \\\hline
		
	\end{tabular}
\end{table}

%%%%%%%%%%%%%%%%%%%%%%%%%%%%%%%%%%%%%%%%%%%%%%%%%%%%%%%%%%%%%%%%%%%%

\section{Proof of Theorem~\ref{CLT:1}}\label{flu}
In this section, we provide a sketch of the proof for Theorem~\ref{CLT:1}. Following a common route for establishing asymptotic normality, we apply Taylor's expansion to the log-likelihood function and try to find the inverse of the corresponding Hessian matrix. However, the randomness of the Hessian matrix introduces additional challenges. The crux of our approach lies in identifying {the corresponding Fisher information matrix} as the weighted Laplacian matrix of appropriate random graphs, which can be carefully analyzed using truncated spectral methods. {Although leveraging the Fisher information matrix to derive asymptotic normality in network models is well established \citep{MR1724040,yan2013central,han2023unified}, existing results require specific model structures, particularly those with logistic-type link functions. To the best of our knowledge, our work provides the first general framework that maintains optimality conditions while accommodating broader model specifications.}

We start by computing the first- and second-order derivatives of the log-likelihood function $l$ for $\u$.
For $i, j\in [n]$ with $i\neq j$, 
\begin{align*}
	\partial_i l(\u) &= \sum_{j\in\delta\{i\}}\frac{\partial_2 f(X_{ij}, u_i-u_j)}{f(X_{ij}, u_i-u_j)}=\sum_{j\in\delta\{i\}} g(X_{ij}, u_i-u_j),\\
	\partial_{ij} l(\u) &= \frac{-\partial_{22}f(X_{ij}, u_i-u_j)\cdot f(X_{ij}, u_i-u_j)+\{\partial_2 f(X_{ij}, u_i-u_j)\}^2}{f^2(X_{ij}, u_i-u_j)}\cdot\mathbb I_{\{(i, j)\in \cE\}}\nonumber\\
	& = -\partial_{2}g(X_{ij}, u_i-u_j)\cdot\mathbb I_{\{(i, j)\in \cE\}}\geq 0\\
	\partial_{ii}l(\u)& = -\sum_{j\in\delta\{i\}}\partial_{ij} l(\u)\leq 0,
\end{align*}
where the first inequality follows from the log-concavity of $f$ that is assumed under the validity assumption. Then the Hessian of $l(\u)$, denoted by $H(\u)$ is given as 
\begin{align*}
	\{H(\u)\}_{i j} &=\left\{\begin{array}{ll}
		\partial_{ij} l(\u) & \quad i\neq j\\
		-\sum_{k\in\delta\{i\}}\partial_{ik} l(\u) & \quad i=j.
	\end{array}\right.
\end{align*}
A crucial observation is that  $H(\u)$ is the negation of an unnormalized graph Laplacian associated with a weighted graph on $[n]$ with edge weights $\partial_{ij} l(\u)$. Conditional on $\cE$, $H(\u)$ is a random quantity due to the comparison outcomes $X_{ij}$. 
As such, we let $H^*(\u) = \E[H(\u)\mid \cE]$. 
It is easy to check that the Fisher information matrix, $-H^*(\u)$ remains an unnormalized weighted graph Laplacian that depends only on $\cG(\cV, \cE)$ and $\u$. 
In the following, we normalize $-H^*(\u^*)$ as
\begin{align}
	-H^*(\u^*) = \D^{1/2}\L_{sym}\D^{1/2} = \D^{1/2}(I- \A)\D^{1/2},\label{lap}
\end{align}
where $\D = \text{diag}\{-H^*(\u^*)\}$ is the diagonal of $-H^*(\u^*)$ (degree matrix), $\L_{sym}$ is the normalized Laplacian of the weighted graph associated with $-H^*(\u^*)$, $\A = I-\L_{sym}$ is the normalized adjacency matrix, and $I$ is the identity matrix. 
By the mean value theorem,
\begin{align*}
	\nabla l(\u^*) &= \nabla l(\u^*) - \nabla l(\uu) = -\bar{H}(\W)(\uu - \u^*)\\
	& = -H^*(\u^*)(\uu - \u^*) + \{H^*(\u^*)-H(\u^*)\}(\uu - \u^*) + \{H(\u^*) - \bar{H}(\W)\}(\uu - \u^*),
\end{align*}
where 
\begin{align*}
	\bar{H}(\W) = \begin{pmatrix}
		\partial_{11}l(\w_1), \ldots, \partial_{1n}l(\w_1)\\
		\\
		\partial_{n1}l(\w_n), \ldots, \partial_{nn}l(\w_n)
	\end{pmatrix}\in\R^{n\times n},
\end{align*} 
with $\W = (\w_1, \ldots, \w_n),$
and $\w_1, \ldots, \w_n$ are $n$ (possibly distinct) points on the line segment between $\uu$ and $\u^*$. 
It is worth noting that $\bar{H}(\W)$ is not a Hessian matrix (it is not even symmetric). 
Heuristically, thanks to the matrix concentration and uniform consistency of $\uu$, $H^*(\u^*)-H(\u^*) = o(H(\u^*))$ and $H(\u^*) - \bar{H}(\W) = o(H(\u^*))$, where the $o(\cdot )$ is interpreted in an appropriate sense that will be made clear in Lemma \ref{mylemma}. 

With some rearrangement and substituting \eqref{lap} into $H^*(\u^*)$, we have 
\begin{align}
	\D^{1/2}(\uu - \u^*) &= -\L^{\dagger}_{sym}\D^{-1/2}\nabla l(\u^*)\label{decomp}\\
	& \ \ \ \ \ +\L^{\dagger}_{sym}\D^{-1/2}\{H^*(\u^*)-H(\u^*)\}(\uu - \u^*)\nonumber\\
	&\ \ \ \ \  + \L^{\dagger}_{sym}\D^{-1/2}\{H(\u^*) - \bar{H}(\W)\}(\uu - \u^*),\nonumber
\end{align}
where $\L_{sym}^\dagger$ is the pseudoinverse of $\L_{sym}$ and can be explicitly expressed using spectral expansion:
\begin{align}
	\L_{sym}^\dagger = \sum_{t=0}^\infty(\A-\mathcal P_1)^t - \mathcal P_1,\label{se}
\end{align}
where $\mathcal P_1$ is the orthogonal projector to the null space of $\L$. Particularly, letting $H_{ij} = \{H^*(\u^*)\}_{ij}$, the entry of $\A$ and $\mathcal P_1$ can be written as 
\begin{align}\label{AP}
	\D_{ii} = \sum_{j\in\delta\{i\}}H_{ij}, \quad \A_{ij} =\frac{H_{ij}}{{ \D_{ii}^{1/2}\D_{jj}^{1/2} }}, \quad \{\mathcal P_1\}_{ij}=\frac{{ \D_{ii}^{1/2}\D_{jj}^{1/2} }}{\sum_{k\in [n]}\D_{kk}}.
\end{align}
Based on the decomposition in \eqref{decomp}, we will proceed to prove the desired asymptotic normality result as follows. 
We will show that the first term on the right-hand side of \eqref{decomp} appropriately normalizes $(\uu - \u^*)$ and gives the asymptotic normality.  
Meanwhile, through a meticulous analysis (Lemma~\ref{mylemma}), it can be established that the remaining terms are of order $o_p(1)$. Consequently, the desired follows by applying Slutsky's theorem. 
\begin{lemma}\label{mylemma}
	Under the assumptions in Theorem \ref{CLT:1}, the following are true:
	\begin{enumerate}
		\item [(1)]$\|\L^{\dagger}_{sym}\D^{-1/2}\{H^*(\u^*)-H(\u^*)\}(\uu - \u^*)\|_\infty = o_p(1)$. 
		\item [(2)]$\|\L^{\dagger}_{sym}\D^{-1/2}\{H(\u^*) - \bar{H}(\W)\}(\uu - \u^*)\|_\infty = o_p(1)$.
		\item [(3)]$\max_{i\in [n]}\sum_{t=2}^\infty \{(\A-\mathcal P_1)^t\}_{ii} = o_p(1)$. 
	\end{enumerate}
\end{lemma}
The proof of Lemma~\ref{mylemma} is rather technical and deferred to the appendix. 
We now use the results in Lemma~\ref{mylemma} to establish the asymptotic normality of $\uu$. 
By statements {\it (1)} and {\it (2)} in Lemma \ref{mylemma} and Slutsky's theorem, it remains to verify for every fixed $S\subset\N$ with $|S| = s$ ($S$ independent of $n$),  $\{\L^{\dagger}_{sym}\D^{-1/2}\nabla l(\u^*)\}_S\to  \mathcal{N}(0, I_s) $.
Since 
\begin{align}
	\L^{\dagger}_{sym}\D^{-1/2}\nabla l(\u^*) = \D^{-1/2}\nabla l(\u^*) + (\L^{\dagger}_{sym}-I)\D^{-1/2}\nabla l(\u^*)\label{12}
\end{align}
and $\{\D^{-1/2}\nabla l(\u^*)\}_S\to \mathcal{N}(0, I_s)$ by the central limit theorem, the proof is finished if $\{(\L^{\dagger}_{sym}-I)\D^{-1/2}\nabla l(\u^*)\}_S = o_p(1)$ componentwise.
To verify this, we apply Chebyshev's inequality. 
Note $\E[(\L^{\dagger}_{sym}-I)\D^{-1/2}\nabla l(\u^*)] = 0$ and the second moment calculation yields that for each $k\in [n]$ and the $k$th canonical basis vector $e_k\in\R^n$, 
\begin{align*}
	&\E\left[e_k^\top(\L^{\dagger}_{sym}-I)\D^{-1/2}\nabla l(\u^*)\{\nabla l(\u^*)\}^\top\D^{-1/2}(\L^{\dagger}_{sym}-I)e_k\right]\\
	=&\  e_k^\top(\L^{\dagger}_{sym}-I)\D^{-1/2}\E[\nabla l(\u^*)\{\nabla l(\u^*)\}^\top]\D^{-1/2}(\L^{\dagger}_{sym}-I)e_k\\
	=&\  -e_k^\top(\L^{\dagger}_{sym}-I)\D^{-1/2}H^*(\u^*)\D^{-1/2}(\L^{\dagger}_{sym}-I)e_k\\
	=&\ e_k^\top(\L^{\dagger}_{sym}-I)\L_{sym}(\L^{\dagger}_{sym}-I)e_k\\
	=&\  (\L^{\dagger}_{sym} + \mathcal P_1-I)_{kk} + (\L_{sym} + \mathcal P_1-I)_{kk} \\
	\stackrel{\eqref{se}}{=}&\  \sum_{t=2}^\infty \{(\A-\mathcal P_1)^t\}_{kk} = o_p(1),
\end{align*} 
where the second equality uses the fact that $\E[\nabla l(\u^*)\{\nabla l(\u^*)\}^\top] = -H^*(\u^*)$ and
the last equality follows from statement {\it (3)} in Lemma \ref{mylemma}. 
The proof is finished by applying Chebyshev's inequality.

\section{Discussion}\label{sec:6}

Despite the broad applicability of the established result in terms of model parameterization, there are practical considerations that motivate further exploration of the problem. For example, in scenarios where multiple data sources are available for the same subjects of interest, different pairwise comparison models might be employed. A promising avenue for future investigation involves integrating the estimated information derived from different yet related models.
The second aspect pertains to temporal dynamics. As can be seen in the ATP example, the data span nearly fifty years. In such contexts, it is more appropriate to also account for the temporal effects on the statistical models employed to gauge comparison outcomes and on the sampling of the graph. Addressing these considerations may involve incorporating additional covariate information.

In addition, one theoretical research direction is to apply our technique for establishing asymptotic normality, that is, the representation of the Fisher information matrix as a weighted graph Laplacian, to a broader class of network models. In particular, we plan to adapt it to the $\beta$-model \citep{chatterjee2011random}, which shares structural similarities with the BT model, and its generalized variants in econometric networks \citep{graham2017econometric,chen2021analysis}. These investigations will form part of our future work.

\appendix
\renewcommand{\theequation}{\Alph{section}.\arabic{equation}} % A.1, B.1, etc.
\setcounter{equation}{0} % Reset equation counter

\section{Appendix: Proof of Lemma \ref{mylemma}}\label{appA}

	We start by making two observations. 
The first states that the spectral gap of $\L_{sym}$ can be bounded using the degree of heterogeneity of the graph. 
The second provides an upper bound for $\D^{-1/2}\{H^*(\u^*)-H(\u^*)\}(\uu - \u^*)$ and $\D^{-1/2}\{H(\u^*) - \bar{H}(\W)\}(\uu - \u^*)$ in the $\ell_\infty$ norm. These two observations are summarized as the following two lemmas respectively. 
Let $\kappa_n = c_{n,4}/c_{n,3}\in [1, \infty)$. 

\begin{lemma}\label{prop1}
	If $\lim_{n\to\infty}(\kappa^{7}_nq^6_n\log n)/(np^7_n)= 0$, then there exists an absolute constant $C>0$ such that
	$\|\A-\mathcal P_1\|_2 \leq 1-C\{p_n/(\kappa_nq_n)\}^2$ holds with probability tending to one.
\end{lemma}

\begin{lemma}\label{prop2}
	If Assumptions \ref{ass:0}--\ref{ass:3} hold, then
	\begin{align*}
		&\|\D^{-1/2}\{H^*(\u^*)-H(\u^*)\}(\uu - \u^*)\|_\infty = O_p\left(\frac{c_{n,2}c_{n,4}^3}{c_{n,3}^{7/2}}\left\{\frac{q_n^{5}(\log n)^4}{np_n^6}\right\}^{1/2}\right)\\
		&\|\D^{-1/2}\{H(\u^*) - \bar{H}(\W)\}(\uu - \u^*)\|_\infty = O_p\left(\frac{c^2_{n,2}c_{n,5}}{c_{n,3}^{5/2}}\left\{\frac{q_n^5(\log n)^6}{np_n^6}\right\}^{1/2}\right) 
	\end{align*}
\end{lemma}

The proofs of Lemmas \ref{prop1}--\ref{prop2} are deferred to the end. 
We now prove Lemma  \ref{mylemma} assuming both Lemmas \ref{prop1}--\ref{prop2} hold.
We begin with the first two statements. 
Let %$L$ be an integer to be specified later and let 
\begin{align*}
	\bm \zeta\in\left\{\D^{-1/2}\{H^*(\u^*)-H(\u^*)\}(\uu - \u^*), \D^{-1/2}\{H(\u^*) - \bar{H}(\W)\}(\uu - \u^*)\right\}.
\end{align*}
Recall that
\begin{align}\label{AP}
	\D_{ii} = \sum_{j\in\delta\{i\}}H_{ij}, \quad \A_{ij} =\frac{H_{ij}}{{\D_{ii}^{1/2} \D_{jj}^{1/2} }}, \quad \{\mathcal P_1\}_{ij}=\frac{ {\D_{ii}^{1/2} \D_{jj}^{1/2} }}{\sum_{k\in [n]}\D_{kk}}.
\end{align}
Note the all-ones vector $\bm 1_n$ is in the zero-eigenspace of $H(\u^*)$ or $H^*(\u^*)$ (in both the left- and right-eigenspace since $H(\u^*)$ and $H^*(\u^*)$ are symmetric). {As a result,}
\begin{align*}
	&\mathcal P_1\D^{-1/2}\{H^*(\u^*)-H(\u^*)\}(\uu - \u^*)\\
	=&\ \frac{1}{\sum_{i\in [n]}\D_{ii}}(\D_{11}^{1/2}, \ldots, \D_{nn}^{1/2})^\top\bm 1^\top_n\{H^*(\u^*)-H(\u^*)\}(\uu - \u^*) = 0.
\end{align*}
On the other hand, since $\nabla l(\u^*) - \nabla l(\uu) = -\bar{H}(\W)(\uu - \u^*)$, 
\begin{align*}
	&-\mathcal P_1\D^{-1/2}\bar{H}(\W)(\uu - \u^*) = \frac{1}{\sum_{i\in [n]}\D_{ii}}(\D_{11}^{1/2}, \ldots, \D_{nn}^{1/2})^\top\bm 1^\top_n(\nabla l(\u^*) - \nabla l(\uu))\\
	=&\ \frac{1}{\sum_{i\in [n]}\D_{ii}}\left\{\sum_{i\in [n]}\sum_{j\in\delta\{i\}}g(X_{ij}, \widehat{u}_i-\widehat{u}_j)-g(X_{ij}, u_i^*-u_j^*)\right\}(\D_{11}^{1/2}, \ldots, \D_{nn}^{1/2})^\top =0,
\end{align*}
which follows as a result of the validity assumption on $f$ (that implies $g(x; y) + g(-x; -y) = 0$). 
Combining both pieces, we obtain $\mathcal P_1\bm\zeta=0$.
Meanwhile, since $\L_{sym}^\dagger = \sum_{t=0}^\infty(\A-\mathcal P_1)^t - \mathcal P_1$, it follows from the Cauchy--Schwarz inequality that
\begin{align}
	\|\L^{\dagger}_{sym}\bm\zeta\|_\infty =\ & \left\|\sum_{t=0}^\infty(\A-\mathcal P_1)^t\bm\zeta\right\|_\infty\\
	\leq\ & \left\|\sum_{t=0}^L(\A-\mathcal P_1)^t\bm\zeta\right\|_\infty+\left\|\sum_{t>L}(\A-\mathcal P_1)^t\bm\zeta\right\|_2\nonumber\\
	\leq\ &  \left\|\sum_{t=0}^L(\A-\mathcal P_1)^t\bm\zeta\right\|_\infty+\sum_{t>L}\|\A-\mathcal P_1\|_2^t\left\|\bm\zeta\right\|_2\nonumber\\
	\leq\ & \left\|\sum_{t=0}^L\A^t\bm\zeta\right\|_\infty+ \frac{{n^{1/2}}\|\A-\mathcal P_1\|_2^L}{1-\|\A-\mathcal P_1\|_2}\left\|\bm\zeta\right\|_{\infty}\nonumber\\
	\leq\ & \left\|\sum_{t=0}^L\A^t\bm\zeta\right\|_\infty+ \frac{{n^{1/2} }\|\A-\mathcal P_1\|_2^L}{1-\|\A-\mathcal P_1\|_2}.\label{lop}
\end{align}
If choosing $L =C(\kappa_nq_n/p_n)^2\log n$ for some sufficiently large constant $C>0$, then the second term in \eqref{lop} is $o_p(1)$ as a result of Lemma \ref{prop1}. 
For the first term in \eqref{lop}, it follows from the direct computation using \eqref{AP} that $\A^t:\R^n\to\R^n$, viewed as an operator in the $\ell_\infty$ norm, can be bounded as 
\begin{align}
	\|\A^t\|_{\ell_\infty\to\ell_\infty} &= \max_{i\in [n]}\sum_{j\in [n]}|\{\A^t\}_{ij}|\label{klopk}\\
	&\stackrel{\eqref{AP}}{\leq}\max_{i \in [n]} \sum_{j_1 \in [n]} \cdots \sum_{j_{t-1} \in [n]}\sum_{j\in [n]} \frac{H_{ij_1}\cdots H_{j_{t-1}j} }{ {\D_{ii}^{1/2}} \D_{j_1j_1}\cdots \D_{j_{t-1}j_{t-1}} {\D_{jj}^{1/2}}}\nonumber\\
	&\leq  \max_{i \in [n]}\frac{1}{\min_{k\in [n]} {\D_{kk}^{1/2} }}\sum_{j_1 \in [n]} \cdots \sum_{j_{t-1} \in [n]} \frac{H_{ij_1}\cdots H_{j_{t-1}j} }{{\D_{ii}^{1/2} } \D_{j_1j_1}\cdots \D_{j_{t-1}j_{t-1}}}\nonumber\\
	& = \frac{\max_{i\in [n]} {\D_{ii}^{1/2} }}{\min_{k\in [n]}{\D_{kk}^{1/2} }}\lesssim\left(\frac{\kappa_nq_n}{p_n}\right)^{1/2},\nonumber
\end{align} 
where the last inequality follows from the standard degree concentration bounds \eqref{ctryb}. Consequently, 
\begin{align*}
	\left\|\sum_{t=0}^L\A^t\bm\zeta\right\|_\infty\lesssim L\left({\frac{\kappa_nq_n}{p_n}}\right)^{1/2}\|\bm\zeta\|_\infty\lesssim (\log n)\cdot\left(\frac{\kappa_n q_n}{p_n}\right)^{5/2}\|\bm\zeta\|_\infty = o_p(1),
\end{align*}
where the last step follows from Lemma \ref{prop2} and the assumptions in the lemma. 

To prove the last statement in Lemma~\ref{mylemma}, note that $\{(\A-\mathcal P_1)^t\}_{ii}\leq\|(\A-\mathcal P_1)^t\|_2$ for every $t$. 
By a similar truncation argument as above, we only need to prove the asymptotic negligibility of the partial sum of $t$ from $2$ to $L$:
\begin{align*}
	\max_{i\in [n]}\sum_{t=2}^L \{(\A-\mathcal P_1)^t\}_{ii} \leq\max_{i\in [n]}\sum_{t=2}^L \{\A^t\}_{ii}\lesssim \frac{L\kappa_{n}^{1/2}}{(np_n)^{1/2} }\lesssim {\kappa_{n}^{5/2}}\left\{\frac{q_n^4(\log n)^2}{np_n^5}\right\}^{1/2} = o_p(1), 
\end{align*}
where the first inequality follows from the fact that $(\A-\mathcal P_1)^t\preceq \A^t$ for $t\geq 1$, and the second inequality follows from similar computation as in \eqref{klopk}, that is,
\begin{align*}
	\{\A^t\}_{ii}
	&{\leq}\max_{i \in [n]} \sum_{j_1 \in [n]} \cdots \sum_{j_{t-1} \in [n]} \frac{H_{ij_1}\cdots H_{j_{t-1}i} }{ {\D_{ii}^{1/2}} \D_{j_1j_1}\cdots \D_{j_{t-1}j_{t-1}} {\D_{ii}^{1/2}}}\nonumber\\
	&\leq  \max_{i \in [n]} (\max_{j_{t-1} \in [n]}  \frac{H_{j_{t-1}i}}{ {\D_{j_{t-1}j_{t-1}} }}) \sum_{j_1 \in [n]} \cdots \sum_{j_{t-1} \in [n]} \frac{H_{ij_1}\cdots H_{j_{t-2}j_{t-1}} }{{\D_{ii} } \D_{j_1j_1}\cdots \D_{j_{t-2}j_{t-2}}}\nonumber\\
	& =   \max_{i \in [n]} (\max_{j_{t-1} \in [n]}  \frac{H_{j_{t-1}i}}{ {\D_{j_{t-1}j_{t-1}} }})  \stackrel{\eqref{ctryb}}{\lesssim} \left(\frac{\kappa_n}{np_n}\right)^{1/2}.\nonumber
\end{align*} 

\section{Proofs of Lemmas \ref{prop1}-\ref{prop2}}
Before showing Lemmas \ref{prop1}-\ref{prop2}, we first cite the uniform consistency result in \cite{han2022general}:
\begin{lemma}\label{thm:cons}
	Let $\cG\sim\G(n, p_n, q_n)$. 
	Under Assumptions \ref{ass:0}-\ref{ass:2}, if
	\begin{align*}
		\alpha_n = \frac{c_{n,2}}{c_{n,3}}\left\{\frac{q_n^2(\log n)^3}{np_n^3}\right\}^{1/2}\to 0\qquad n\to\infty,  
	\end{align*}
	then there exists $C>0$ (independent of $n$) such that, with probability one, the MLE $\widehat{\u}$ uniquely exists and $\|\widehat{\u}-\u^*\|_\infty\leq C\alpha_n$  for all sufficiently large $n$. 
\end{lemma}

\subsection{Proof of Lemma \ref{prop1}}\label{ssss}
The proof of Lemma~\ref{prop1} follows a similar approach as \cite[Lemma 5.3]{han2023unified}. 
Denote the eigenvalues of $\mathcal{L}_{\text {sym }} = I-\mathcal A = I-\mathcal{D}^{-1 / 2} \mathcal{W} \mathcal{D}^{-1 / 2}$ in increasing order as $0=\lambda_1 \leq \cdots \leq \lambda_n \leq 2$, where $\D$ and $\mathcal{W}$ are all random matrices (where randomness is from the comparison graph $\cG$ only). In particular, $\mathcal{W}$ is an off-diagonal weight matrix with $\mathcal{W}_{ij}=\mathbb{I}_{\{(i, j)\in \cE\}} z_{ij}$, where 
\begin{align*}
	z_{i j}&= \E_{X_{ij}}\left[-\partial_{2}g(X_{ij}, u^*_i-u^*_j)\right]\in [c_{n,3}, c_{n, 4}], \quad i\neq j
\end{align*}
due to Assumption \ref{ass:2}, with $\E_{X_{ij}}[\cdot ]$ denoting expectation for $X_{ij}$. It is clear $\left\|\mathcal{A}-\mathcal{P}_1\right\|_2=\max \left\{1-\lambda_2,  \lambda_n-1\right\}$.

To analyze the spectrum of $\mathcal{L}_{\text {sym }}$, we consider its expected version first and then apply concentration inequalities.
Let $\overline{\mathcal{D}}$ and $\overline{\mathcal{W}}$ be the expectation of $\mathcal{D}$ and $\mathcal{W}$, respectively, and $\bar{\mathcal{L}}_{\text {sym }}=I-\bar{\mathcal{D}}^{-1 / 2} \bar{\mathcal{W}} \bar{\mathcal{D}}^{-1 / 2}$. 
By definition, $\bar{\mathcal{L}}_{\text {sym }}$ is also a normalized graph Laplacian matrix; we denote the eigenvalues of $\overline{\mathcal{L}}_{\text {sym }}$ by $0=\bar{\lambda}_1 \leq \cdots \leq \bar{\lambda}_n \leq 2$. 
A computation based on Cheeger's inequalities yields a lower bound for the spectral gap of the expected Laplacian \citep{MR1421568, bauer_jost_2013}: There exists an absolute constant $c>0$ such that
\begin{align}
	\max \left\{1-\bar{\lambda}_2, \bar{\lambda}_n-1\right\}\leq 1-c\left(\frac{p_n}{\kappa_n q_n}\right)^2. \label{esg}
\end{align}
We next apply a concentration argument to show that the spectra of $\mathcal{L}_{\text {sym }}$ and $\bar{\mathcal{L}}_{\text {sym }}$ are close. 
To this end, we first apply Weyl's inequality and the triangle inequality to obtain 
\begin{align*}
	&\max\left\{ \left| (1 - \lambda_2) - (1 - \bar{\lambda}_2) \right|, \ \left| (\lambda_n - 1) - (\bar{\lambda}_n - 1) \right| \right\} \\
	\leq\ &\left\| \mathcal{L}_{\text{sym}} - \overline{\mathcal{L}}_{\text{sym}} \right\|_2 
	= \left\| \mathcal{D}^{-1/2} \mathcal{W} \mathcal{D}^{-1/2} - \overline{\mathcal{D}}^{-1/2} \overline{\mathcal{W}} \overline{\mathcal{D}}^{-1/2} \right\|_2 \\
	\leq\ &\left\| \mathcal{D}^{-1/2} - \overline{\mathcal{D}}^{-1/2} \right\|_2 \left\| \mathcal{W} \right\|_2 \left\| \mathcal{D}^{-1/2} \right\|_2 + \left\| \overline{\mathcal{D}}^{-1/2} \right\|_2 \left\| \mathcal{W} - \overline{\mathcal{W}} \right\|_2 \left\| \mathcal{D}^{-1/2} \right\|_2 \\
	&+ \left\| \overline{\mathcal{D}}^{-1/2} \right\|_2 \left\| \overline{\mathcal{W}} \right\|_2 \left\| \mathcal{D}^{-1/2} - \overline{\mathcal{D}}^{-1/2} \right\|_2
\end{align*}
By the Chernoff bound, with probability at least $1-n^{-2}$, all $D_{i i}$ are concentrated around their means
\begin{align}
	\left|\mathcal{D}_{i i}-\overline{\mathcal{D}}_{i i}\right| \lesssim \left({c_{n,4}\overline{\mathcal{D}}_{i i} \log n}\right)^{1/2} \quad i \in[n],\label{ctryb}
\end{align}
so that
\begin{align*}
	\left\|\mathcal{D}^{-1/2}-\overline{\mathcal{D}}^{-1/2}\right\|_2 & =\max _{i \in[n]}\left|\mathcal{D}_{i i}^{-1/2}-\overline{\mathcal{D}}_{i i}^{-1/2}\right| \lesssim  \max _{i \in[n]}\left(\overline{\mathcal{D}}_{i i}^{-3/2}\left|\mathcal{D}_{i i}-\overline{\mathcal{D}}_{i i}\right|\right) \lesssim  \left\{\frac{c_{n,4}\log n}{(c_{n,3}np_n)^2}\right\}^{1/2} \\
	\left\|\mathcal{D}^{-1/2}\right\|_2 & \lesssim \left\|\overline{\mathcal{D}}^{-1/2}\right\|_2 \lesssim  \sqrt{\frac{1}{c_{n,3}np_n}}.
\end{align*}
On the other hand,
\begin{align*}
	\|\overline{\mathcal{W}}\|_2=\left\|\overline{\mathcal{D}}^{1/2}  (I - \bar{\mathcal{L}}_{\text {sym}})\overline{\mathcal{D}}^{1/2}\right\|_2 \leq\|\overline{\mathcal{D}}\|_2\left\|I - \bar{\mathcal{L}}_{\text {sym}}\right\|_2 \leq \|\overline{\mathcal{D}}\|_2= \max _{i\in [n]} \overline{\mathcal{D}}_{i i} \lesssim c_{n,4}nq_n.
\end{align*}
By the matrix Bernstein inequality \citep{tropp2012user}, with probability at least $1-n^{-2}$,
\begin{align*}
	\|\mathcal{W}-\overline{\mathcal{W}}\|_2 \lesssim c_{n,4}\left({nq_n \log n}\right)^{1/2}.
\end{align*}
In conclusion, the upper bound for $\left\|\mathcal{L}_{\text {sym}}-\overline{\mathcal{L}}_{\text {sym}}\right\|_2 $ is 
\begin{align}
	\left\|\mathcal{L}_{\text {sym }}-\overline{\mathcal{L}}_{\text {sym}}\right\|_2 \lesssim \left(\frac{\kappa_n^3 q_n^2\log n}{np_n^3}\right)^{1/2} + \left(\frac{\kappa_n^2 q_n\log n}{np_n^2}\right)^{1/2}.\label{flc}
\end{align}
The desired result follows if the expected spectral gap in \eqref{esg} dominates the concentration deviations \eqref{flc}:	\begin{align*}
	\frac{ \left(\frac{\kappa_n^3 q_n^2\log n}{np_n^3}\right)^{1/2} + \left(\frac{\kappa_n^2 q_n\log n}{np_n^2}\right)^{1/2}}{\left(\frac{p_n}{\kappa_n q_n}\right)^2}\to 0\quad n\to\infty,
\end{align*}
which holds true if $\lim_{n\to \infty}(\kappa_n^7q_n^6\log n)/(np_n^7) \to 0$. 

\subsection{Proof of Lemma \ref{prop2}}
Recall the formula for the Hessian matrix $H(\u)$ and $\bar{H}(\W)$:
\begin{align}\label{tww}
	\{H(\u)\}_{i j} &=\left\{\begin{array}{ll}
		\partial_{ij} l(\u) = -\partial_{2}g(X_{ij}, u_i-u_j)\times \mathbb I_{\{(i, j)\in \cE\}} & \quad i\neq j\\
		-\sum_{k\in\delta\{i\}}\partial_{ik} l(\u) & \quad i=j
	\end{array}\right.\\
	\{\bar{H}(\W)\}_{i j} &=\left\{\begin{array}{ll}
		\partial_{ij} l(\bm w_i) = -\partial_{2}g(X_{ij}, w_{ii}-w_{ij})\times \mathbb I_{\{(i, j)\in \cE\}} & \quad i\neq j\\
		-\sum_{k\in\delta\{i\}}\partial_{ik} l(\bm w_i) & \quad i=j
	\end{array}\right.
\end{align}
where $\w_1, \ldots, \w_n$ are $n$ points on the line segment between $\uu$ and $\u^*$, and $\w_i = (w_{i1}, \ldots, w_{in})\in\R^n$.   
Write the $i$th element of $\{\bar{H}(\W)-H\left(\u^*\right)\}\left(\widehat{\u}-\u^*\right)$ as
\begin{align*}
	&\left|\left[\{\bar{H}(\W)-H\left(\u^*\right)\}\left(\widehat{\u}-\u^*\right)\right]_i\right| \\
	=&\left|\sum_{j\neq i}\left\{\bar{H}(\W)-H\left(\u^*\right)\right\}_{ij}(\widehat{u}_j-u^*_j)+\left\{\bar{H}(\W)-H\left(\u^*\right)\right\}_{ii}(\widehat{u}_i-u^*_i)\right|\\
	=&\left|\sum_{j\neq i}\left\{\bar{H}(\W)-H\left(\u^*\right)\right\}_{ij}(\widehat{u}_j-u^*_j)-\sum_{j\neq i}\left\{\bar{H}(\W)-H\left(\u^*\right)\right\}_{ij}(\widehat{u}_i-u^*_i)\right|  \\
	=&\left|\sum_{j \neq i}\left\{\bar{H}(\W)-H\left(\u^*\right)\right\}_{ij}\left\{\left(\widehat{u}_j-u_j^*\right)-\left(\widehat{u}_{i}-u_{i}^*\right)\right\}\right|
\end{align*}
By the mean value theorem, for every $i, j$ and all sufficiently large $n$, there exists $\xi_{ij}\in [-M_n-1, M_n+1]$ depending on $i, j$ such that
\begin{align*}
	\left|\left\{\bar{H}(\W)-H\left(\u^*\right)\right\}_{ij}\right|&\stackrel{\eqref{tww}}{=} \left|\left\{\partial_{2}g(X_{ij}, u^*_i-u^*_j)-\partial_{2}g(X_{ij}, w_{ii}-w_{ij})\right\}\times \mathbb I_{\{(i, j)\in \cE\}}\right|\\
	& = \left|\partial_{22}g(X_{ij}, \xi_{ij})\left\{\left(u_i^*-w_{ii}\right)-\left(u_j^*-w_{ij}\right)\right\}\times \mathbb I_{\{(i, j)\in \cE\}}\right|\\
	&\leq 2c_{n,5}\|\bm w_i-\u^*\|_{\infty}\times \mathbb I_{\{(i, j)\in \cE\}}\\
	&\leq 2c_{n,5}\|\uu-\u^*\|_{\infty}\times \mathbb I_{\{(i, j)\in \cE\}},
\end{align*}
where the penultimate inequality follows from Assumption \ref{ass:3}.
Summing over $j\in [n]$, 
\begin{align*}
	\left|\left[\left\{\bar{H}(\W)-H\left(\u^*\right)\right\}\left(\widehat{\u}-\u^*\right)\right]_i\right|
	\leq 4c_{n,5}\|\widehat{\u}-\u^*\|^2_{\infty}|\delta\{i\}|.
\end{align*}
Combining the results with Lemma \ref{thm:cons}, we can have
\begin{align*}
	&\left|\left[D^{-1/2}\left\{\bar{H}(\W)-H\left(\u^*\right)\right\}\left(\widehat{\u}-\u^*\right)\right]_i\right|\leq[D^{-1/2}]_{ii}\left|\left[\left\{\bar{H}(\W)-H\left(\u^*\right)\right\}\left(\widehat{\u}-\u^*\right)\right]_i\right|\\
	&\leq  4c_{n,5}[D^{-1/2}]_{ii}\|\widehat{\u}-\u^*\|^2_{\infty} \left({|\delta\{i\}|}\right)^{1/2} \leq \frac{c^2_{n,2}c_{n,5}}{c^{5/2}_{n,3}}\frac{q_n^{5/2}(\log n)^3}{n^{1/2}p_n^3}.
\end{align*}
The proof of the first statement is completed.
Next, we prove the second statement in Lemma \ref{prop2}. 
Similar to the previous proof, we can decompose $\left \lvert[\{H(\u^*) - H^\ast(\u^\ast)\}(\uu - \u^*)]_i\right \lvert $ into two parts and obtain
\begin{align*}
	&\left \lvert[\{H(\u^*) - H^\ast(\u^\ast)\}(\uu - \u^*)]_i\right \lvert\\
	\leq&\left|\sum_{j\neq i}\left\{H(\u^*) - H^\ast(\u^\ast)\right\}_{ij}(\widehat{u}_j-u^*_j)\right|+\left|\sum_{j\neq i}\left\{H(\u^*) - H^\ast(\u^\ast)\right\}_{ij}(\widehat{u}_i-u^*_i)\right|\\
	\leq&\underbrace{\left|\sum_{j\neq i}\left\{H(\u^*) - H^\ast(\u^\ast)\right\}_{ij}(\widehat{u}_j-u^*_j)\right|}_{(A.1)}+\underbrace{\left|\sum_{j\neq i}\left\{H(\u^*) - H^\ast(\u^\ast)\right\}_{ij}\right|\|\uu-\u^*\|_{\infty}}_{(A.2)}.
\end{align*}
To bound $(A.2)$, note $\sum_{j\neq i}\left\{H(\u^*) - H^\ast(\u^\ast)\right\}_{ij}$ is a sum of independent centered random variables bounded by $2c_{n,4}$. By Hoeffding's inequality, with probability at least $1-n^{-2}$,  
\begin{align*}
	\left|\sum_{j\neq i}\left\{H(\u^*) - H^\ast(\u^\ast)\right\}_{ij}\right|\lesssim c_{n,4}  \left(|\delta\{i\}| \log n\right)^{1/2} 	\end{align*}
Therefore, 
\begin{align}\label{A2bound}
	(A.2)\lesssim c_{n,4} \left(|\delta\{i\}| \log n\right)^{1/2} 
	\times \|\uu-\u^\ast\|_{\infty}\lesssim  \left({|\delta\{i\}|}\right)^{1/2}   \times \frac{c_{n,2}c_{n,4} q_n(\log n)^2}{c_{n,3}np_n^{3/2}}. 
\end{align}
Unfortunately, $(A.1)$ cannot be treated similarly as a result of dependence between \\$\left\{H(\u^*) - H^\ast(\u^\ast)\right\}_{ij}$ and $\widehat{u}_j-u^*_j$. 
To address this issue, we apply a leave-one-out approach based on a modification of \cite{gao2023uncertainty}. 
For $i \in [n]$, let
\begin{align}
	l^{(-i)}(\u) = \sum_{j\neq i, k\neq i, (j, k)\in \cE}\log f(X_{jk},u_j-u_k).
\end{align}
We define $\widehat{\u}^{(-i)}\in\R^{n-1}$ as leave-one-out estimator by maximizing the likelihood 
\begin{align*}
	\uu^{(-i)} =\argmax_{\u\in \R^{n-1}: \bm 1_{n-1}^\top \u = -u^*_i}l^{(-i)}(\u).
\end{align*}
We also denote $\u_{-i}\in\R^{n-1}$ as $\u$ deleting the $i$th component.   
Note that $\uu^{(-i)}$ is also uniformly consistent, that is,
\begin{align*}
	\|\uu^{(-i)}-\u_{-i}^\ast\|_{\infty} \lesssim \frac{c_{n,2}}{c_{n,3}} \left\{\frac{q_n^2(\log n)^3}{np_n^3}\right\}^{1/2}.
\end{align*}
Now the term $(A.1)$ can be bounded using the triangle inequality as
\begin{align}\label{A1bound}
	(A.1)\leq &\underbrace{\left|\sum_{j\neq i}\left\{H(\u^*) - H^\ast(\u^\ast)\right\}_{ij}\left(\widehat u_j-\widehat u_j^{(i)}\right)\right|}_{ (J.1)}+\underbrace{\left|\sum_{j\neq i}\left\{H(\u^*) - H^\ast(\u^\ast)\right\}_{ij}\left(\widehat u^{(i)}_j- u_j^\ast\right)\right|}_{(J.2)}.
\end{align}
Since $\uu^{(-i)}$ is independent of $\{X_{ij}\}_{j\in\delta\{i\}}$, we can bound ($J.2$) similar to the steps in bounding $(A.2)$ by first conditioning on $\uu^{(-i)}$.
With probability at least $1-n^{-2}$, 
\begin{align}\label{J2bound}
	(J.2) \lesssim c_{n,4}  \left(|\delta\{i\}| \log n\right)^{1/2}  %\sqrt{\log n\cdot|\delta\{i\}|}
	\times \|\uu^{(-i)}-\u_{-i}^\ast\|_{\infty}.
\end{align}
To bound ($J.1$), we apply the Cauchy--Schwarz inequality to obtain
\begin{align}
	(J.1)&\leq \left\|\left[\left\{H(\u^*) - H^\ast(\u^\ast)\right\}_{i,:}\right]_{-i}\right\|_2\|\uu_{-i}-\uu^{(-i)}\|_2\nonumber\\
	&\lesssim c_{n,4} \left({|\delta\{i\}|}\right)^{1/2} \times \|\uu_{-i}-\uu^{(-i)}\|_2,\label{J1bound}
\end{align}
where $\left\{H(\u^*) - H^\ast(\u^\ast)\right\}_{i,:}$ denotes the $i$th row of $H(\u^*) - H^\ast(\u^\ast)$. 
To further upper bound $\|\uu_{-i}-\uu^{(-i)}\|_2$, we use a Taylor expansion as follows: 
\begin{align*}
	(\widehat{\u}^{(-i)} - \widehat{\u}_{-i})^\top\nabla l^{(-i)}(\widehat{\u}_{-i}) &=  -(\widehat{\u}^{(-i)} - \widehat{\u}_{-i})^\top\left\{\nabla l^{(-i)}(\widehat{\u}^{(-i)}) - \nabla l^{(-i)}(\widehat{\u}_{-i})\right\}\\
	&= (\widehat{\u}^{(-i)} - \widehat{\u}_{-i})^\top  \{-\nabla^2 l^{(-i)}(\widetilde{\u})\} (\widehat{\u}^{(-i)} - \widehat{\u}_{-i}),
\end{align*}
where $\widetilde{\u}$ lies on the line segment between $\widehat{\u}^{(-i)}$ and $\widehat{\u}_{-i}$.  
Applying the Cauchy--Schwarz inequality, 
\begin{align}\label{total}
	\|\widehat{\u}^{(-i)} - \widehat{\u}_{-i}\|_2 \leq \frac{\lVert \nabla l^{(-i)}(\widehat{\u}_{-i})\lVert_2}{\lambda_{2} \{-\nabla^2 l^{(-i)}(\widetilde{\u})\} },
\end{align}
where $\lambda_{2}\{\cdot\}$ (with some abuse of notation) denotes the second smallest eigenvalue.
Since $ \nabla l(\widehat{\u}) = 0$, separating the terms involving $i$ and the rest apart, we have
\begin{align*}
	&\|\nabla l^{(-i)}(\widehat{\u}_{-i})\|_2^2 = \sum_{j\neq i, (i, j)\in \cE} g(X_{ij},\widehat{u}_i-\widehat{u}_j)^2\\
	&\lesssim \sum_{j\neq i, (i, j)\in \cE} \left\{g(X_{ij},\widehat{u}_i-\widehat{u}_j) - g(X_{ij},u^*_i-u^*_j) \right\}^2 + \sum_{j\neq i, (i, j)\in \cE} g(X_{ij},u^*_i-u^*_j)^2\\
	&\lesssim c_{n,4} nq_n \|\uu-\u^\ast\|_{\infty}^2 + \sum_{j\neq i, (i, j)\in \cE} g(X_{ij},u^*_i-u^*_j)^2.
\end{align*} 
Under Assumption \ref{ass:1}, $\{g(X_{ij},u^*_i-u^*_j)\}_{j: (i, j)\in \cE}$ are subgaussian random variables with subgaussian norms uniformly bounded by $c_{n,2}$. 
According to \cite[Lemma 2.7.6]{vershynin2018high}, $\{g(X_{ij},u^*_i-u^*_j)^2\}_{j: (i, j)\in \cE}$ are sub-exponential random variables with sub-exponential norms uniformly bounded by $c_{n,2}^2$. 
It follows from \cite[Proposition 2.7.1, Exercise 2.7.10]{vershynin2018high} that 
\begin{align*}
	&\max_{j: (i, j)\in \cE}\E[g(X_{ij},u^*_i-u^*_j)^2]\lesssim c_{n,2}^2\\
	&\max_{j: (i, j)\in \cE}\|g(X_{ij},u^*_i-u^*_j)^2-\E[g(X_{ij},u^*_i-u^*_j)^2] \|_{\psi_1}\lesssim c_{n,2}^2,
\end{align*}
where $\|\cdot\|_{\psi_1}$ denotes the sub-exponential norm. 
By Bernstein's inequality, 
\begin{align*}
	&\sum_{j\neq i, (i, j)\in \cE} \{g(X_{ij},u^*_i-u^*_j) \}^2\\
	& \lesssim \sum_{j\neq i, (i, j)\in \cE} \E[g(X_{ij},u^*_i-u^*_j)^2] \\
	& \qquad \qquad + \left|\sum_{j\neq i, (i, j)\in \cE} \left[g(X_{ij},u^*_i-u^*_j)^2-\E[g(X_{ij},u^*_i-u^*_j)^2]\right]\right|\\
	&\leq 2\sum_{j\neq i, (i, j)\in \cE} \E[g(X_{ij},u^*_i-u^*_j)^2]\lesssim c_{n,2}^2 nq_n
\end{align*}
holds with probability at least $1-n^{-3}$. Consequently, with probability at least $1-n^{-3}$, 
\begin{align}\label{num}
	\|\nabla l^{(-i)}(\widehat{\u}_{-i})\|_2 &\lesssim \left\{(c_{n,4} \|\uu-\u^\ast\|_{\infty}^2 + c_{n,2}^2)nq_n\right\}^{1/2}  \lesssim c_{n,2}(nq_n)^{1/2},
\end{align}
where the last step follows from the observation 
\begin{align*}
	c_{n,4} \|\uu-\u^\ast\|_{\infty}^2\lesssim\frac{c_{n,2}^2c_{n,4}}{c_{n,3}^2}\frac{q_n^2\log n}{np_n^3}\to 0,
\end{align*}
as $n \to \infty.$

On the other hand, by a similar calculation as in the proof of Lemma \ref{prop1} (for normalized Laplacian), with probability at least $1-n^{-2}$, 
\begin{align}\label{dem}
	\lambda_{2} \{\nabla^2 l^{(-i)}(\widetilde{\u})\} \gtrsim c_{n,3}np_n\left(\frac{p_n}{\kappa_n q_n}\right)^2,
\end{align}
where the extra term $c_{n,3}np_n$ stems from a lower bound for the singular value of the normalization matrix. 
Combining \eqref{total}--\eqref{dem}, we conclude
\begin{align}
	\|\widehat{\u}^{(-i)} - \widehat{\u}_{-i}\|_2 \lesssim \frac{c_{n,2}c_{n,4}^2 q_n^{5/2}}{c_{n,3}^3p_n^3n^{1/2}}.\label{leaveoneout}
\end{align}
Putting \eqref{A1bound}, \eqref{J2bound}, \eqref{J1bound}, and \eqref{leaveoneout} together, 
\begin{align}
	(A.1)&\lesssim c_{n,4} \left(|\delta\{i\}| \log n\right)^{1/2} %\sqrt{\log n\cdot|\delta\{i\}|}
	\left(\|\uu^{(-i)}-\u_{-i}^\ast\|_{\infty}+\|\widehat{\u}^{(-i)} - \widehat{\u}_{-i}\|_2\right)\nonumber\\
	&\lesssim  \left({|\delta\{i\}|}\right)^{1/2} \times \frac{c_{n,2}c_{n,4}^3 q_n^{5/2}(\log n)^2}{c_{n,3}^3p_n^3n^{1/2}}.\label{kaizi}
\end{align}
Consequently, combining \eqref{A2bound} and \eqref{kaizi}, we have 
\begin{align*}
	|D^{-1/2}\{H(\u^*) - H^\ast(\u^\ast)\}(\uu - \u^*)|_i &= D_{ii}^{-1/2}\left \lvert[\{H(\u^*) - H^\ast(\u^\ast)\}(\uu - \u^*)]_i\right \lvert\\
	&\leq D_{ii}^{-1/2}[(A.1)+(A.2)]\\
	&\lesssim \frac{c_{n,2}c_{n,4}^3 q_n^{5/2}(\log n)^2}{c_{n,3}^{7/2}p_n^3n^{1/2}}.
\end{align*}
The second part of Lemma \ref{prop2} is proved.

\section{Proofs of Corollary \ref{myco} and Theorems \ref{individual_error}-\ref{CI}}

\subsection{Proof of Corollary \ref{myco}}

\begin{proof}
	If $M^*<\infty$, then $\sup_{n}c_{n,1}<1$ so that Assumption~\ref{ass:0} holds under condition $(\log n)^8/ (np_n)\to 0$. 
	Since $\bbA$ is finite and $g(x;y)$ is continuous in $y$, $g(x;y)$ is uniformly bounded over $\bbA\times [-M^*, M^*]$, implying $\sup_nc_{n,2}<\infty$. 
	Additionally, by the strict log-concavity of $f$, $\partial_2 g(x; y)>0$, which combined with the continuity of $\partial_2 g(x; y)$, the boundedness of $M^*$, and the finiteness of $\bbA$ implies $0<\inf_n c_{n,3}\leq \sup_n c_{n,4}<\infty$. 
	A similar argument shows $\sup_n c_{n,5}<\infty$. 
	The proof is finished by appealing to \eqref{good}. 
\end{proof}

\subsection{Proof of Theorem \ref{individual_error}}
According to \eqref{decomp} and \eqref{12}, we have
\begin{align*}
	\D^{1/2}(\uu - \u^*) &=  -\D^{-1/2}\nabla l(\u^*) - (\L^{\dagger}_{sym}-I)\D^{-1/2}\nabla l(\u^*)\\
	& +\L^{\dagger}_{sym}\D^{-1/2}\{H^*(\u^*)-H(\u^*)\}(\uu - \u^*)\nonumber\\
	&+ \L^{\dagger}_{sym}\D^{-1/2}\{H(\u^*) - \bar{H}(\W)\}(\uu - \u^*).\nonumber
\end{align*}
Then, for any $k \in [n]$, we have 
\begin{align*}
	\{\D^{1/2}(\uu - \u^*)\}_k &=  -\{\D^{-1/2}\nabla l(\u^*)\}_k + o_p(1),
\end{align*}
via Lemma \ref{mylemma} and Chebyshev's inequality. Since $\D$ is diagonal matrix,
\begin{align*}
	\widehat{u}_k - u^*_k = - \D^{-1}_{kk}\{\nabla l(\u^*)\}_k + o_p(\D^{-1/2}_{kk}).
\end{align*}
By Hoeffding's inequality, we obtain
\begin{align*}
	\D^{-1}_{kk} = O_p\left( \{c_{n,3}|\delta\{i\}|\}^{-1}\right), \  \{\nabla l(\u^*)\}_k = O_p(c_{n,2}\{ |\delta\{i\}| \cdot \log n\}^{1/2}).
\end{align*}
Therefore,
\begin{align*}
	|\widehat{u}_k - u^*_k| = O_p\left(\frac{c_{n,2}}{c_{n,3}}\left\{\frac{\log n}{|\delta\{i\}|}\right\}^{1/2} \right).
\end{align*}
We finish the proof of Theorem  \ref{individual_error}.	
\subsection{Proof of Theorem \ref{CI}}
According to Slutsky's theorem and Theorem \ref{CLT:1}, it suffices to show that for any fixed $i \in [n]$ that
\begin{align*}
	\frac{\{\rho_{i}(\u^*)\}^{-1} - \{\rho_{i}(\uu)\}^{-1}}{ \{\rho_{i}(\u^*)\}^{-1}} = o_p(1).
\end{align*}
Note that $g(x,y) =  \partial_2(\log f(x; y))$, we have
\begin{align*}
	\{\rho_{i}(\u^*)\}^{-1} = \sum_{j\in\delta\{i\}}\int_{\bbA} \{g(x; {u}^*_{i}-{u}^*_{j})\}^2 f(x; {u}^*_{i}-{u}^*_{j}) \d x
\end{align*}
and
\begin{eqnarray*}
	&&\{\rho_{i}(\u^*)\}^{-1} - \{\rho_{i}(\uu)\}^{-1} \\
	&=& \sum_{j\in\delta\{i\}}\int_{\bbA} \{g(x; {u}^*_{i}-{u}^*_{j})\}^2 f(x; {u}^*_{i}-{u}^*_{j}) - \{g(x; \widehat{u}_{i}-\widehat{u}_{j})\}^2 f(x; \widehat{u}_{i}-\widehat{u}_{j}) \d x.
\end{eqnarray*}
For simplicity, let $ g_{ij}^*(x) := g(x; {u}^*_{i}-{u}^*_{j})$, $f_{ij}^*(x) := f(x; {u}^*_{i}-{u}^*_{j})$,  $\widehat{g}_{ij}(x) := g(x; \widehat{u}_{i}-\widehat{u}_{j})$ and $\widehat{f}_{ij}(x) := f(x; \widehat{u}_{i}-\widehat{u}_{j})$.
It is enough to show that for any $j \in [n]$, 
\begin{align*}
	\frac{\int_{\bbA} \{g_{ij}^*(x) \}^2 f_{ij}^*(x) - \{\widehat{g}_{ij}(x)\}^2 \widehat{f}_{ij}(x)  \d x}{ \int_{\bbA} \{g_{ij}^*(x)\}^2 f_{ij}^*(x)\d x} = o_p(1).
\end{align*}	
According to H\"older's inequality,
\begin{align*}
	\left| \frac{\int_{\bbA} \{g_{ij}^*(x) \}^2 f_{ij}^*(x) - \{\widehat{g}_{ij}(x)\}^2 \widehat{f}_{ij}(x)  \d x}{ \int_{\bbA} \{g_{ij}^*(x)\}^2 f_{ij}^*(x) \d x} \right| \leq \max_{x \in \bbA}\left| \frac{ \{g_{ij}^*(x) \}^2 f_{ij}^*(x) - \{\widehat{g}_{ij}(x)\}^2 \widehat{f}_{ij}(x) }{  \{g_{ij}^*(x)\}^2 f_{ij}^*(x)} \right|.
\end{align*}
For the right-hand side,
\begin{align*}
	\frac{ \{g_{ij}^*(x) \}^2 f_{ij}^*(x) - \{\widehat{g}_{ij}(x)\}^2 \widehat{f}_{ij}(x) }{  \{g_{ij}^*(x)\}^2 f_{ij}^*(x)} = \underbrace{ \frac{ \{g_{ij}^*(x) \}^2 - \{\widehat{g}_{ij}(x)\}^2}{  \{g_{ij}^*(x)\}^2}}_{\Gamma_1}  + \underbrace{  \frac{\{\widehat{g}_{ij}(x)\}^2}{\{g_{ij}^*(x)\}^2 } \cdot \frac{  f_{ij}^*(x) -  \widehat{f}_{ij}(x) }{  f_{ij}^*(x)}}_{\Gamma_2}.
\end{align*}
On the one hand, we decompose the term $\Gamma_1$ into two parts,
\begin{align*}
	\Gamma_1 =  \frac{ g_{ij}^*(x) - \widehat{g}_{ij}(x) }{  g_{ij}^*(x)} \cdot \frac{ g_{ij}^*(x) + \widehat{g}_{ij}(x) }{ g_{ij}^*(x)}.
\end{align*}
By the mean value theorem,
\begin{align*}
	\left| \frac{ g_{ij}^*(x) - \widehat{g}_{ij}(x) }{  g_{ij}^*(x)} \right|=  \left| \frac{\partial_2 g(x; w_{ij})}{ g_{ij}^*(x)}  \right|  \cdot \left| ({u}^*_{i}-{u}^*_{j}) -  (\widehat{u}_{i}-\widehat{u}_{j}) \right|,
\end{align*}
where $w_{ij}$	 is the intermediate value between $({u}^*_{i}-{u}^*_{j}) $ and $ (\widehat{u}_{i}-\widehat{u}_{j})$. Therefore, according to Lemma \ref{thm:cons}, we have with probability approaching one that
\begin{align*}
	\left| \frac{ g_{ij}^*(x) - \widehat{g}_{ij}(x) }{  g_{ij}^*(x)} \right| \lesssim \frac{c_{n,2}c_{n,4}}{c_{n,3}}\left\{\frac{q_n^2(\log n)^3}{np_n^3}\right\}^{1/2}=o(1),
\end{align*}	
where the last step follows assuming \eqref{good} holds. Meanwhile, 
\begin{align*}
	\left| \frac{ g_{ij}^*(x) + \widehat{g}_{ij}(x) }{ g_{ij}^*(x)}\right| = \left|2+ \frac{ \widehat{g}_{ij}(x) -g_{ij}^*(x)}{ g_{ij}^*(x)}\right| = O_p(1).
\end{align*}
Therefore, $| \Gamma_1| = o_p(1).$ On the other hand, $ \{\widehat{g}_{ij}(x)/g_{ij}^*(x)\}^2 = O_p(1)$ following a similar approach as above. By the same reasoning, 
\begin{align*}
	\log (f_{ij}^*(x)) - \log (\widehat{f}_{ij}(x)) =   g(x; v_{ij})\{ ({u}^*_{i}-{u}^*_{j}) -  (\widehat{u}_{i}-\widehat{u}_{j}) \}.
\end{align*}
Then, based on Assumption \ref{ass:1}, Lemma \ref{thm:cons}, and Hoeffding's inequality, the following inequlity
\begin{align*}
	|\log (f_{ij}^*(x)) - \log (\widehat{f}_{ij}(x)) | \leq  \frac{c_{n,2}^2}{c_{n,3}}\left\{\frac{q_n^2(\log n)^4}{np_n^3}\right\}^{1/2} = o(1).
\end{align*}
holds with probability approaching one. As a result,
\begin{align*}
	\left|\frac{  f_{ij}^*(x) -  \widehat{f}_{ij}(x) }{  f_{ij}^*(x)} \right|= o_p(1).
\end{align*}
Therefore, $| \Gamma_2| = o_p(1).$ The above proof is shown for fixed $x$. However, it is straightforward to extend it to the uniform case since we only rely on the event that Lemma \ref{thm:cons} holds. Consequently,
\begin{align*}
	\max_{x \in \bbA}\left| \frac{ \{g_{ij}^*(x) \}^2 f_{ij}^*(x) - \{\widehat{g}_{ij}(x)\}^2 \widehat{f}_{ij}(x) }{  \{g_{ij}^*(x)\}^2 f_{ij}^*(x)} \right|	= o_p(1).\end{align*} 
We finish the proof of Theorem \ref{CI}.

\section*{Funding}
 R. Han was supported by the Hong Kong Research Grants Council (No. 14301821) and the Hong Kong Polytechnic University (P0044617, P0045351). 
 
 \section*{Relevant code}
 The code to reproduce the main results in this paper is available at \url{https://github.com/RJ-HAN-STAT/Ranking_Code}.
 
\bibliographystyle{apalike}
\bibliography{paper-ref}

@article{mcfadden1973conditional,
  title={Conditional logit analysis of qualitative choice behavior},
  author={McFadden, D},
  journal={Frontiers in Economics},
  pages={105--142},
  year={1973},
  publisher={Academic Press}
}

@article{chen2021analysis,
	title={Analysis of networks via the sparse $\beta$-model},
	author={Chen, Mingli and Kato, Kengo and Leng, Chenlei},
	journal={J. Roy. Statist. Soc. Ser. B},
	volume={83},
	number={5},
	pages={887--910},
	year={2021},
	publisher={Oxford University Press}
}

@article{wapman2022quantifying,
	title={Quantifying hierarchy and dynamics in {US} faculty hiring and retention},
	author={Wapman, K Hunter and Zhang, Sam and Clauset, Aaron and Larremore, Daniel B},
	journal={Nature},
	volume={610},
	number={7930},
	pages={120--127},
	year={2022},
	publisher={Nature Publishing Group UK London}
}

@article{10.1111/rssa.12124,
	author = {Varin, Cristiano and Cattelan, Manuela and Firth, David},
	title = {Statistical Modelling of Citation Exchange Between Statistics Journals},
	journal={J. Roy. Statist. Soc. Ser. A},
	volume = {179},
	number = {1},
	pages = {1-63},
	year = {2015},
	month = {11},
	issn = {0964-1998},
	doi = {10.1111/rssa.12124},
	url = {https://doi.org/10.1111/rssa.12124},
	eprint = {https://academic.oup.com/jrsssa/article-pdf/179/1/1/49348123/jrsssa\_179\_1\_1.pdf},
}

@article{10.1257/pandp.20221063,
	Author = {Gu, Jiaying and Koenker, Roger},
	Title = {Ranking and Selection from Pairwise Comparisons: Empirical Bayes Methods for Citation Analysis},
	Journal = {AEA Papers and Proceedings},
	Volume = {112},
	Year = {2022},
	Month = {May},
	Pages = {624–29},
	DOI = {10.1257/pandp.20221063},
	URL = {https://www.aeaweb.org/articles?id=10.1257/pandp.20221063}}

@inproceedings{Christiano2017,
  title={Deep reinforcement learning from human preferences},
  author={Christiano, Paul F. and Leike, Jan and Brown, Tom and Martic, Miljan and Legg, Shane and Amodei, Dario},
  booktitle={NeurIPS},
  pages={4299--4307},
  year={2017}
}

@inproceedings{Rafailov2023,
  title={Direct Preference Optimization: Your Language Model is Secretly a Reward Model},
  author={Rafailov, Rafael and Sharma, Archit and Mitchell, Eric and Manning, Christopher D. and Ermon, Stefano and Finn, Chelsea},
  booktitle={NeurIPS},
  year={2023},
  url={https://arxiv.org/abs/2305.18290}
}

@inproceedings{Sun2025,
  title={Rethinking Reward Modeling in Preference-based Large Language Model Alignment},
  author={Sun, Hao and Shen, Yunyi and Ton, Jean-Francois},
  booktitle={ICLR},
  year={2025},
  url={https://openreview.net/forum?id=rfdblE10qm}
}

@article{collingwood2022evaluating,
  title={Evaluating the effectiveness of different player rating systems in predicting the results of professional snooker matches},
  author={Collingwood, James AP and Wright, Michael and Brooks, Roger J},
  journal={Eur. J. Oper. Res.},
  volume={296},
  number={3},
  pages={1025--1035},
  year={2022},
  publisher={Elsevier}
}

@article{baker2014dynamic,
  title={A dynamic paired comparisons model: Who is the greatest tennis player?},
  author={Baker, Rose D and McHale, Ian G},
	journal={Eur. J. Oper. Res.},
  volume={236},
  number={2},
  pages={677--684},
  year={2014},
  publisher={Elsevier}
}

@article{angelini2022weighted,
	title={Weighted Elo rating for tennis match predictions},
	author={Angelini, Giovanni and Candila, Vincenzo and De Angelis, Luca},
	journal={Eur. J. Oper. Res.},
	volume={297},
	number={1},
	pages={120--132},
	year={2022},
	publisher={Elsevier}
}

@article{eco1958,
	ISSN = {00129682, 14680262},
	URL = {http://www.jstor.org/stable/1907622},
	author = {Gerard Debreu},
	journal = {Econometrica},
	number = {3},
	pages = {440--444},
	publisher = {[Wiley, Econometric Society]},
	title = {Stochastic Choice and Cardinal Utility},
	urldate = {2025-05-11},
	volume = {26},
	year = {1958}
}

@article{lentz2023anatomy,
	title={The Anatomy of Sorting—Evidence From Danish Data},
	author={Lentz, Rasmus and Piyapromdee, Suphanit and Robin, Jean-Marc},
	journal={Econometrica},
	volume={91},
	number={6},
	pages={2409--2455},
	year={2023},
	publisher={Wiley Online Library}
}

@article{bozoki2016application,
	title={An application of incomplete pairwise comparison matrices for ranking top tennis players},
	author={Boz{\'o}ki, S{\'a}ndor and Csat{\'o}, L{\'a}szl{\'o} and Temesi, J{\'o}zsef},
	journal={Eur. J. Oper. Res.},
	volume={248},
	number={1},
	pages={211--218},
	year={2016},
	publisher={Elsevier}
}

@article{newman2023efficient,
  title={Efficient computation of rankings from pairwise comparisons},
  author={Newman, MEJ},
  journal={J. Mach. Learn. Res.},
  volume={24},
  number={238},
  pages={1--25},
  year={2023}
}

@article{qu2023sinkhorn,
  title={On Sinkhorn's Algorithm and Choice Modeling},
  author={Qu, Zhaonan and Galichon, Alfred and Ugander, Johan},
  journal={arXiv preprint arXiv:2310.00260},
  year={2023}
}

@article{chen2023note,
	title={A note on statistical inference for noisy incomplete 1-bit matrix},
	author={Chen, Yunxiao and Li, Chengcheng and Ouyang, Jing and Xu, Gongjun},
	journal={J. Mach. Learn. Res.},
	volume={24},
	pages={1--66},
	year={2023}
}

@inproceedings{hendrickx2020minimax,
  title={Minimax rate for learning from pairwise comparisons in the BTL model},
  author={Hendrickx, Julien and Olshevsky, Alex and Saligrama, Venkatesh},
  booktitle={ICML},
  pages={4193--4202},
  year={2020},
  organization={PMLR}
}

@inproceedings{chen2015spectral,
  title={Spectral mle: Top-k rank aggregation from pairwise comparisons},
  author={Chen, Yuxin and Suh, Changho},
  booktitle={ICML},
  pages={371--380},
  year={2015},
  organization={PMLR}
}

@inproceedings{bong2022generalized,
  title={Generalized results for the existence and consistency of the MLE in the {B}radley--{T}erry--{L}uce model},
  author={Bong, Heejong and Rinaldo, Alessandro},
  booktitle={ICML},
  pages={2160--2177},
  year={2022},
  organization={PMLR}
}

@article{han2023unified,
  title={A unified analysis of likelihood-based estimators in the {P}lackett--{L}uce model},
  author={Han, Ruijian and Xu, Yiming},
  journal={Ann. Stat.},
  volume = {53},
  number = {5},
  pages = {2077--2102},
  year={2025}
}

@article{loewen2012testing,
  title={Testing the power of arguments in referendums: A {B}radley--{T}erry approach},
  author={Loewen, Peter John and Rubenson, Daniel and Spirling, Arthur},
  journal={Elect. Stud.},
  volume={31},
  number={1},
  pages={212--221},
  year={2012},
  publisher={Elsevier}
}

@article{tropp2012user,
	    AUTHOR = {Tropp, Joel A.},
  title={User-friendly tail bounds for sums of random matrices},
  JOURNAL = {Found. Comput. Math.},
FJOURNAL = {Foundations of Computational Mathematics. The Journal of the
Society for the Foundations of Computational Mathematics},
VOLUME = {12},
YEAR = {2012},
NUMBER = {4},
PAGES = {389--434},
MRCLASS = {60B20 (60F10 60G42 60G50)},
MRNUMBER = {2946459}
}

@article{bauer_jost_2013,
	AUTHOR = {Bauer, Frank and Jost, J\"{u}rgen},
	TITLE = {Bipartite and neighborhood graphs and the spectrum of the
	normalized graph {L}aplace operator},
	JOURNAL = {Comm. Anal. Geom.},
	FJOURNAL = {Communications in Analysis and Geometry},
	VOLUME = {21},
	YEAR = {2013},
	NUMBER = {4},
	PAGES = {787--845},
	MRCLASS = {05C50 (05C81 31C20 58C40)},
	MRNUMBER = {3078942},
	MRREVIEWER = {Yufei Huang},
}

@article{chen2022optimal,
	AUTHOR = {Chen, Pinhan and Gao, Chao and Zhang, Anderson Y.},
	TITLE = {Optimal full ranking from pairwise comparisons},
	JOURNAL = {Ann. Statist.},
	FJOURNAL = {The Annals of Statistics},
	VOLUME = {50},
	YEAR = {2022},
	NUMBER = {3},
	PAGES = {1775--1805},
	MRCLASS = {62F07},
	MRNUMBER = {4441140},

}

@article{gao2023uncertainty,
	AUTHOR = {Gao, Chao and Shen, Yandi and Zhang, Anderson Y.},
	TITLE = {Uncertainty quantification in the {B}radley--{T}erry--{L}uce
	model},
	JOURNAL = {Inf. Inference},
	FJOURNAL = {Information and Inference. A Journal of the IMA},
	VOLUME = {12},
	YEAR = {2023},
	NUMBER = {2},
	PAGES = {1073--1140},
	MRCLASS = {62F07 (62F12 62F25 62H12)},
	MRNUMBER = {4565758},
}

@article{thurstone1927method,
  title={The method of paired comparisons for social values.},
  author={Thurstone, Louis L},
  journal = {J. Abnorm. Psychol.},
  FJOURNAL={The Journal of Abnormal and Social Psychology},
  volume={21},
  number={4},
  pages={384},
  year={1927},
  publisher={American Psychological Association}
}

@article{han2022general,
  title={A general pairwise comparison model for extremely sparse networks},
  author={Han, Ruijian and Xu, Yiming and Chen, Kani},
   JOURNAL = {J. Amer. Statist. Assoc.},
FJOURNAL = {Journal of the American Statistical Association},
  pages={2422--2432},
  	VOLUME = {118},
  NUMBER = {544},
  year={2023},
  publisher={Taylor \& Francis}
}

@article {MR0070925,
	AUTHOR = {Bradley, Ralph Allan and Terry, Milton E.},
	TITLE = {Rank analysis of incomplete block designs. {I}. {T}he method
	of paired comparisons},
	JOURNAL = {Biometrika},
	FJOURNAL = {Biometrika},
	VOLUME = {39},
	YEAR = {1952},
	PAGES = {324--345},
	MRCLASS = {62.0X},
	MRNUMBER = {0070925},
	MRREVIEWER = {W. S. Connor},
}

@article {MR217963,
	AUTHOR = {Rao, P. V. and Kupper, L. L.},
	TITLE = {Ties in paired-comparison experiments: {A} generalization of
	the {B}radley-{T}erry model},
	JOURNAL = {J. Amer. Statist. Assoc.},
	FJOURNAL = {Journal of the American Statistical Association},
	VOLUME = {62},
	YEAR = {1967},
	PAGES = {194--204},
	MRCLASS = {62.64},
	MRNUMBER = {217963},
	MRREVIEWER = {W. A. Thompson, Jr.},
}

@article{davidson1970extending,
	title={On extending the {B}radley--{T}erry model to accommodate ties in paired comparison experiments},
	author={Davidson, Roger R},
	JOURNAL = {J. Amer. Statist. Assoc.},
	volume={65},
	number={329},
	pages={317--328},
	year={1970},
	publisher={Taylor \& Francis Group}
}

@article{MR1724040,
    AUTHOR = {Simons, Gordon and Yao, Yi-Ching},
TITLE = {Asymptotics when the number of parameters tends to infinity in
the {B}radley-{T}erry model for paired comparisons},
JOURNAL = {Ann. Statist.},
FJOURNAL = {The Annals of Statistics},
  volume={27},
  number={3},
  pages={1041--1060},
  year={1999},
  publisher={Institute of Mathematical Statistics},
  MRCLASS = {62J15 (62E20 62F05 62F12)},
  MRNUMBER = {1724040}
}

@article{MR2987494,
    AUTHOR = {Yan, Ting and Yang, Yaning and Xu, Jinfeng},
TITLE = {Sparse paired comparisons in the {B}radley-{T}erry model},
JOURNAL = {Statist. Sinica},
FJOURNAL = {Statistica Sinica},
VOLUME = {22},
YEAR = {2012},
NUMBER = {3},
PAGES = {1305--1318},
MRCLASS = {62J15},
MRNUMBER = {2987494}
}

@article{MR3613103,
    AUTHOR = {Negahban, Sahand and Oh, Sewoong and Shah, Devavrat},
TITLE = {Rank centrality: ranking from pairwise comparisons},
JOURNAL = {Oper. Res.},
FJOURNAL = {Operations Research},
VOLUME = {65},
YEAR = {2017},
NUMBER = {1},
PAGES = {266--287},
MRCLASS = {91B14 (91B06)},
MRNUMBER = {3613103}
}

@article {MR3953449,
	AUTHOR = {Chen, Yuxin and Fan, Jianqing and Ma, Cong and Wang, Kaizheng},
	TITLE = {Spectral method and regularized {MLE} are both optimal for
	top-{$K$} ranking},
	JOURNAL = {Ann. Statist.},
	FJOURNAL = {The Annals of Statistics},
	VOLUME = {47},
	YEAR = {2019},
	NUMBER = {4},
	PAGES = {2204--2235},
	MRCLASS = {62F07 (60J20 62B10)},
	MRNUMBER = {3953449}
}

@article{han2020asymptotic,
    AUTHOR = {Han, Ruijian and Ye, Rougang and Tan, Chunxi and Chen, Kani},
TITLE = {Asymptotic theory of sparse {B}radley--{T}erry model},
JOURNAL = {Ann. Appl. Probab.},
FJOURNAL = {The Annals of Applied Probability},
VOLUME = {30},
YEAR = {2020},
NUMBER = {5},
PAGES = {2491--2515},
MRCLASS = {60F05 (62E20 62F12 62J15)},
MRNUMBER = {4149535}
}

@article {MR3012434,
	AUTHOR = {Cattelan, Manuela},
	TITLE = {Models for paired comparison data: a review with emphasis on
	dependent data},
	JOURNAL = {Statist. Sci.},
	FJOURNAL = {Statistical Science. A Review Journal of the Institute of
	Mathematical Statistics},
	VOLUME = {27},
	YEAR = {2012},
	NUMBER = {3},
	PAGES = {412--433},
	MRCLASS = {62J15},
	MRNUMBER = {3012434},
}

@article{agresti1992analysis,
	title={Analysis of ordinal paired comparison data},
	author={Agresti, Alan},
	journal={J. Roy. Statist. Soc. Ser. C},
	volume={41},
	number={2},
	pages={287--297},
	year={1992},
	publisher={Wiley Online Library}
}

@article{MR3504618,
      AUTHOR = {Shah, Nihar B. and Balakrishnan, Sivaraman and Bradley, Joseph
  and Parekh, Abhay and Ramchandran, Kannan and Wainwright,
  Martin J.},
  TITLE = {Estimation from pairwise comparisons: sharp minimax bounds
  with topology dependence},
  JOURNAL = {J. Mach. Learn. Res.},
  FJOURNAL = {Journal of Machine Learning Research (JMLR)},
  VOLUME = {17},
  YEAR = {2016},
  PAGES = {Paper No. 58, 47},
  MRCLASS = {62J15 (62C20 62F07)},
  MRNUMBER = {3504618}

}

@article {MR2051012,
	AUTHOR = {Hunter, David R.},
	TITLE = {M{M} algorithms for generalized {B}radley--{T}erry models},
	JOURNAL = {Ann. Statist.},
	FJOURNAL = {The Annals of Statistics},
	VOLUME = {32},
	YEAR = {2004},
	NUMBER = {1},
	PAGES = {384--406},
	MRCLASS = {62F07 (62C20 65C20)},
	MRNUMBER = {2051012}
}

@article {MR1545015,
	AUTHOR = {Zermelo, E.},
	TITLE = {Die {B}erechnung der {T}urnier-{E}rgebnisse als ein
	{M}aximumproblem der {W}ahrscheinlichkeitsrechnung},
	JOURNAL = {Math. Z.},
	FJOURNAL = {Mathematische Zeitschrift},
	VOLUME = {29},
	YEAR = {1929},
	NUMBER = {1},
	PAGES = {436--460},
	MRCLASS = {DML},
	MRNUMBER = {1545015},
}

@article {MR0125031,
	AUTHOR = {Erd\H{o}s, Paul and R\'{e}nyi, Alfr\'{e}d},
	TITLE = {On the evolution of random graphs},
	JOURNAL = {Magyar Tud. Akad. Mat. Kutat\'{o} Int. K\"{o}zl.},
	FJOURNAL = {A Magyar Tudom\'{a}nyos Akad\'{e}mia. Matematikai Kutat\'{o} Int\'{e}zet\'{e}nek
	K\"{o}zlem\'{e}nyei},
	VOLUME = {5},
	YEAR = {1960},
	PAGES = {17--61},
	MRCLASS = {05.40},
	MRNUMBER = {0125031},
	MRREVIEWER = {John Riordan},
}

@book{MR1421568,
  title={Spectral graph theory},
  author={Chung, Fan RK},
  volume={92},
  year={1997},
  publisher={American Mathematical Soc.}
}

@book{vershynin2018high,
  title={High-dimensional probability: An introduction with applications in data science},
  author={Vershynin, Roman},
  volume={47},
  year={2018},
  publisher={Cambridge University Press}
}

@article{mosteller2006remarks,
	title={Remarks on the method of paired comparisons},
	author={Mosteller, Frederick},
	journal={Psychometrika},
	volume={16},
	number={3},
	pages={279--289},
	year={1951},
	publisher={Springer}
}

@article{chatterjee2011random,
	title={Random graphs with a given degree sequence},
	author={Chatterjee, Sourav and Diaconis, Persi and Sly, Allan and others},
	journal={The Annals of Applied Probability},
	volume={21},
	number={4},
	pages={1400--1435},
	year={2011},
	publisher={Institute of Mathematical Statistics}
}

@article{graham2017econometric,
	title={An econometric model of network formation with degree heterogeneity},
	author={Graham, Bryan S},
	journal={Econometrica},
	volume={85},
	number={4},
	pages={1033--1063},
	year={2017},
	publisher={Wiley Online Library}
}

@article{yan2013central,
	title={A central limit theorem in the $\beta$-model for undirected random graphs with a diverging number of vertices},
	author={Yan, Ting and Xu, Jinfeng},
	journal={Biometrika},
	volume={100},
	number={2},
	pages={519--524},
	year={2013},
	publisher={Oxford University Press}
}

@article{MR4560195,
	AUTHOR = {Liu, Yue and Fang, Ethan X. and Lu, Junwei},
	TITLE = {Lagrangian inference for ranking problems},
	JOURNAL = {Oper. Res.},
	FJOURNAL = {Operations Research},
	VOLUME = {71},
	YEAR = {2023},
	NUMBER = {1},
	PAGES = {202--223},
	MRCLASS = {62F07 (62B10 62C20)},
	MRNUMBER = {4560195}
}

@article{doi:10.1080/01621459.2024.2316364,
	AUTHOR = {Fan, Jianqing and Lou, Zhipeng and Wang, Weichen and Yu,
	Mengxin},
	TITLE = {Ranking inferences based on the top choice of
	multiway comparisons},
	JOURNAL = {J. Amer. Statist. Assoc.},
	FJOURNAL = {Journal of the American Statistical Association},
	VOLUME = {120},
	YEAR = {2025},
	NUMBER = {549},
	PAGES = {237--250},
	ISSN = {0162-1459},
	MRNUMBER = {4893554}
}

\end{document}